\documentclass{article}
\usepackage{amssymb,amsmath,bm,geometry,graphics,graphicx,url,color,hyperref}
\def\no{\noindent}
\def\pmatrix{\left(\begin{array}}
\def\endpmatrix{\end{array}\right)}

\def\CC{\mathbb{C}}

\def\RR{\mathbb{R}}

\def\D{{\cal D}}

\def\H{{\cal H}}
\def\I{{\cal I}}

\def\P{{\cal P}}

\def\dd{\mathrm{d}}

\def\diag{\mathrm{diag}}

\def\ii{\mathrm{i}}
\def\Re{\mathrm{Re}}
\def\Im{\mathrm{Im}}

\def\sech{\mathrm{sech}}
\def\vec{\mathrm{vec}}

\newtheorem{theo}{Theorem}
\newtheorem{lem}{Lemma}
\newtheorem{cor}{Corollary}
\newtheorem{rem}{Remark}

\def\proof{\noindent\underline{Proof}\quad}
\def\QED{\mbox{~$\Box{~}$}}

\def\bff{{\bm{f}}}
\def\bfg{{\bm{g}}}

\def\bfp{{\bm{p}}}
\def\bfq{{\bm{q}}}

\def\bfw{{\bm{w}}}

\def\bfy{{\bm{y}}}

\def\bfzero{{\bm{0}}}

\def\bfdelta{{\bm{\delta}}}
\def\bfeta{{\bm{\eta}}}

\def\tD{\tilde{D}}
\def\eps{\varepsilon}

\def\mH{\mathrm{H}}

\begin{document}

\title{Spectrally accurate space-time solution of Manakov systems}

\author{Luigi Barletti$\,^{\rm a}$ \and Luigi Brugnano$\,^{\rm a,*}$  \and  Yifa Tang$\,^{\rm b}$ \and Beibei Zhu$\,^{\rm c}$}

\date{\small a) Dipartimento di Matematica e Informatica ``U.\,Dini'', Universit\`a di Firenze, 50134 Firenze, Italy.\\ 
b) Academy of Mathematics and Systems Science, Chinese Academy of Sciences, 100190 Beijing, China.\\
c) National Center for Mathematics and Interdisciplinary Sciences, Academy of Mathematics and Systems Science, Chinese Academy of Sciences, Beijing 100190, China.\\
*) Corresponding author: \url{luigi.brugnano@unifi.it}}

\maketitle

\begin{abstract}  In this paper, we study the numerical solution of Manakov systems by using a spectrally accurate Fourier decomposition in space, coupled with a spectrally accurate time integration. This latter relies on the use of spectral Hamiltonian Boundary Value Methods. The used approach allows to conserve all the physical invariants of the systems. Some  numerical tests are reported.

\medskip
\no{\bf Keywords:~ }Hamiltonian PDEs, Manakov system, energy-conserving methods, Hamiltonian Boundary Value Methods, HBVMs, spectral methods.

\medskip
\no{\bf MSC:~}  65P10, 65M70, 65M20, 65L05, 65L06.

\end{abstract}

\section{Introduction}\label{intro} Generally speaking, the Manakov system (which is a particular case of the vector nonlinear Schr\"o\-din\-ger equation) describes situations where two or more wave modes, of different frequency and/or polarization, are nonlinearly coupled. 
This applies to electromagnetic waves (e.g.\ in optical fibers \cite{Agrawal,MW1997,MAS2012,MEA2013}) as well as to matter waves (e.g.\ in Bose-Einstein condensates \cite{Frantz2010}). In its basic form, the Manakov system is known to be completely integrable \cite{ZM76} and admits vector $N$-soliton solutions \cite{RL95}.

Because of its mathematical relevance and its importance for applications, Manakov systems have been the subject of both theoretical  (see references above and also \cite{CZ2012,Frisquetal2015,RW2013,ZW2017}) and numerical investigation \cite{CC2017,FS2015}. Also the development of novel numerical methods for their solution has been considered by some authors. As an example, Galerkin methods have been considered in \cite{Ismail2008}, compact space-time schemes have been studied in \cite{KHJZ2015}, multi-sympelctic schemes are given in \cite{QSC2014}, and a fourth-order in space -- second-order in time energy-conserving scheme has been defined in \cite{KWHZW2019}. The latter two methods have the noticeable feature of conserving discrete counterparts of some physical conservation laws associated with the systems (see Theorem~\ref{cons} below) but, on the other hand, they have a quite low order in time (second-order at most). Therefore, looking for arbitrarily high-order numerical methods in both space and time, allowing to conserve such invariants, motivates this paper.

A rather general form of the one-dimensional Manakov system is
\begin{equation}\label{ms}
\ii\partial_t \psi = -\beta \partial_{xx}\psi -\gamma|\psi|^2\circ \psi, \qquad (x,t)\in\D \equiv [a,b]\times[0,T],
\end{equation}
where $\psi:\D\rightarrow\CC^n$, $\circ$ is the Hadamard (i.e., componentwise) product, and
\begin{equation}\label{gbp}
\gamma = \left(\gamma_{jk}\right) = \gamma^\top, ~ \beta = \diag(\,\beta_1,\,\dots\,,\beta_n\,)~ \in~\RR^{n\times n},\qquad
|\psi|^2 = \pmatrix{ccc} |\psi_1|^2, &\dots~,&|\psi_n|^2\endpmatrix^\top,
\end{equation} 
being $\psi_j$ the $j$-th component of $\psi$. As a notational convention, hereafter for any given vector $v$, we shall denote 
\begin{equation}\label{v2}
v^2=v\circ v.
\end{equation}
In order to better exemplify our notations, let us write down explicitly the $j$-th component of Eq.\,\eqref{ms}:
$$
 \ii\partial_t \psi_j = -\beta_j \partial_{xx}\psi_j -\sum_{k=1}^n \gamma_{jk}|\psi_k|^2\, \psi_j, \qquad j=1,\dots,n.
$$
Equation (\ref{ms}) can be also written as
\begin{equation}\label{Hop0}
\ii\partial_t \psi = \mathrm{H}(\psi)\psi, 
\end{equation}
with
\begin{equation}\label{Hop}
\mH(\psi) = -\beta \partial_{xx} -\gamma|\psi|^2\circ
\end{equation}
the associated ($|\psi|^2$-dependent) Hamiltonian operator.
We also introduce the Hermitian products $\langle\cdot,\cdot\rangle$ and $\langle\cdot,\cdot\rangle_1$,
\begin{equation}\label{Hprods}
\langle\psi,\phi\rangle = \sum_{j=1}^n \langle\psi_j,\phi_j\rangle_1 \equiv \sum_{j=1}^n \int_a^b \psi_j(x)\overline\phi_j(x)\dd x,
\end{equation}
with the associated norms, 
\begin{equation}\label{norms}
\|\psi_j\|_1^2 = \langle\psi_j,\psi_j\rangle_1, \qquad \|\psi\|^2 = \langle\psi,\psi\rangle \equiv \sum_{j=1}^n\|\psi_j\|_1^2.
\end{equation}
Equation (\ref{ms}) is completed with the initial conditions
\begin{equation}\label{ms1}
\psi(x,0) = \psi^0(x), \qquad x\in[a,b],
\end{equation}
and periodic boundary conditions. We shall assume that $\psi^0(x)$ is enough regular (as a periodic function), so as to ensure a corresponding suitably regular (periodic in space) solution. We next show that problem (\ref{ms}) and (\ref{ms1}) has a number of conserved quantities, during the evolution.

\begin{theo}\label{cons} With reference to (\ref{Hop})--(\ref{norms}), the following quantities are conserved for the solution of (\ref{ms}) and (\ref{ms1}):
\begin{description}
\item[mass of the $j$-th component:]
\begin{equation}\label{Mj}
M_j(t) = \|\psi_j\|_1^2 \equiv \int_a^b |\psi_j(x,t)|^2\dd x,\qquad j=1,\dots,n,
\end{equation}
\item[total mass:]
\begin{equation}\label{M}
M(t) = \|\psi\|^2 \equiv \sum_{j=1}^n M_j(t),
\end{equation}
\item[total momentum:]
\begin{equation}\label{K}
K(t) = \Im\langle\partial_x \psi,\psi\rangle \equiv \frac{1}{2\ii}\sum_{j=1}^n \int_a^b \left[\overline\psi_j(x,t)\partial_x\psi_j(x,t)-\psi_j(x,t) \partial_x\overline\psi_j(x,t)\right]\dd x,
\end{equation}
\item[energy:] 
\begin{eqnarray}\nonumber
E(t) &=& -\frac{1}2\langle \beta\partial_{xx}\psi,\psi\rangle - \frac{1}4\langle\gamma|\psi|^2\circ\psi,\psi\rangle 
\equiv \frac{1}2\langle \beta\partial_x\psi,\partial_x\psi\rangle - \frac{1}4\langle\gamma|\psi|^2\circ\psi,\psi\rangle\\ \label{E}
&\equiv& \frac{1}2\sum_{j=1}^n\int_a^b\left[
\beta_j|\partial_x\psi_j(x,t)|^2-\frac{1}2\sum_{k=1}^n\gamma_{jk}|\psi_k(x,t)|^2|\psi_j(x,t)|^2\right]\dd x.
\end{eqnarray}
\end{description}
 \end{theo}
 \proof
 By using the $j$-th component of system (\ref{ms}), i.e.,\footnote{Hereafter, when not necessary, we shall omit the arguments of the functions, for sake of brevity.}
 $$\ii\partial_t \psi_j = \mH_j(\psi)\psi_j,$$
 we obtain\,\footnote{Hereafter, the $\dot{~}$ will denote the total time derivative.}
 \begin{eqnarray*}
 \dot M_j &=& \langle\partial_t \psi_j,\psi_j\rangle_1 + \langle\psi_j,\partial_t \psi_j\rangle_1 = -\ii\langle\mH_j(\psi)\psi_j,\psi_j\rangle_1+\ii\langle\psi_j,\mH_j(\psi)\psi_j\rangle_1\\[1mm] &=&-\ii\langle\mH_j(\psi)\psi_j,\psi_j\rangle_1+\ii\langle\mH_j(\psi)\psi_j,\psi_j\rangle_1=0,
 \end{eqnarray*}
 due to the hermitianity of $\mH_j(\psi)$ with periodic boundary conditions, w.r.t. the inner product $\langle\cdot,\cdot\rangle_1$, thus proving the conservation of (\ref{Mj}). Consequently,  (\ref{M}) is conserved as well.
 
 Concerning (\ref{K}), one has:
 \begin{eqnarray*}
 \dot K &=& \Im\left( \langle\partial_x\partial_t\psi,\psi\rangle+\langle\partial_x\psi,\partial_t\psi\rangle\right)
 = \Im\left( -\langle\partial_t\psi,\partial_x\psi\rangle+\langle\partial_x\psi,\partial_t\psi\rangle\right)\\
 &=& \Im\left( -\overline{\langle\partial_x\psi,\partial_t\psi\rangle}+\langle\partial_x\psi,\partial_t\psi\rangle\right)
 = 2\Im\left(\langle\partial_x\psi,\partial_t\psi\rangle\right)\\
  &=& 2\Im\left(\ii \langle\partial_x\psi,\mH(\psi)\psi\rangle\right) = 2\Re\left(\langle\partial_x\psi,\mH(\psi)\psi\rangle\right).
 \end{eqnarray*}
 Considering that
 $$
 \langle\partial_x\psi,\mH(\psi)\psi\rangle = -\langle\partial_x\psi,\beta\partial_{xx}\psi\rangle-\langle\partial_x\psi,\gamma|\psi|^2\circ\psi\rangle,
 $$
 and
$$
 \langle\partial_x\psi,\beta\partial_{xx}\psi\rangle = -\langle\beta\partial_{xx}\psi,\partial_x\psi\rangle = -\overline{\langle\partial_x\psi,\beta\partial_{xx}\psi\rangle}, 
 $$
 and also, using the symmetry of $\gamma$,
 \begin{eqnarray*}
 \langle\partial_x\psi,\gamma|\psi|^2\circ\psi\rangle &=& \sum_{j,k=1}^n\gamma_{jk}\int_a^b\partial_x\psi_j\overline\psi_j|\psi_k|^2\dd x
 = - \sum_{j,k=1}^n\gamma_{jk}\int_a^b\psi_j\partial_x\overline\psi_j|\psi_k|^2\dd x\\
 && - \sum_{j,k=1}^n\gamma_{jk}\int_a^b|\psi_j|^2\partial_x|\psi_k|^2\dd x = - \sum_{j,k=1}^n\gamma_{jk}\int_a^b\psi_j\partial_x\overline\psi_j|\psi_k|^2\dd x\\
 &&\underbrace{-\frac{1}2 \sum_{j,k=1}^n\gamma_{jk}\int_a^b\partial_x\left(|\psi_j|^2|\psi_k|^2\right)\dd x}_{=0} = 
 -\overline{\langle\partial_x\psi,\gamma|\psi|^2\circ\psi\rangle},
 \end{eqnarray*}
 one has that $\langle\partial_x\psi,\beta\partial_{xx}\psi\rangle$  and  $\langle\partial_x\psi,\gamma|\psi|^2\circ\psi\rangle$, and therefore 
 $ \langle\partial_x\psi,\mH(\psi)\psi\rangle$, are imaginary quantities, from which $\dot K=0$ follows.
 
 At last, for the conservation of (\ref{E}), one has:
 $$2\dot E = -\langle\beta\partial_{xx}\partial_t\psi,\psi\rangle-\langle\beta\partial_{xx}\psi,\partial_t\psi\rangle -\frac{1}2\sum_{j,k=1}^n\gamma_{jk}
 \int_a^b \partial_t\left(\psi_k\overline\psi_k\psi_j\overline\psi_j\right)\dd x.
 $$ 
 By using the symmetry of $\gamma$, we can write
 $$
 \frac{1}2\sum_{j,k=1}^n\gamma_{jk} \int_a^b \partial_t\left(\psi_k\overline\psi_k\psi_j\overline\psi_j\right)\dd x = 
 \sum_{j,k=1}^n  \gamma_{jk}\int_a^b |\psi_k|^2\overline\psi_j\partial_t\psi_j\dd x+\sum_{j,k=1}^n  \gamma_{jk}\int_a^b |\psi_k|^2\psi_j\partial_t\overline\psi_j\dd x.$$
Consequently, using also the hermitianity of $\beta \partial_{xx}$ and (\ref{Hop}), one has:
 \begin{eqnarray*}
 2\dot E&=&-\langle\partial_t\psi,\beta\partial_{xx}\psi\rangle-\overline{\langle\partial_t\psi,\beta\partial_{xx}\psi\rangle}
 -\langle\partial_t\psi,\gamma|\psi|^2\circ\psi\rangle-\overline{\langle\partial_t\psi,\gamma|\psi|^2\circ\psi\rangle}\\[2mm]
 &=&\langle\partial_t\psi,\mathrm{H}(\psi)\psi\rangle+\overline{\langle\partial_t\psi,\mathrm{H}(\psi)\psi\rangle}.
 \end{eqnarray*}
 Substituting, see (\ref{Hop0}), $\partial_t\psi$ by $-\ii\mathrm{H}(\psi)\psi$ we immediately arrive at $\dot E=0$.\,\QED
\bigskip

In this paper, after a suitable space semi-discretization, we shall derive energy-conserving methods for numerically solving problem (\ref{ms}) and (\ref{ms1}), when periodic boundary conditions are prescribed. Moreover, the proposed space-time method has the potentiality of providing  {\em spectral accuracy} in both space and time. With these premises, the structure of the paper is as follows: in Section~\ref{rform} we transform the problem (\ref{ms}) and (\ref{ms1}) into real form, then provide a space expansion along a suitable functional basis; in Section~\ref{hbvms} we describe the time integration through arbitrarily high-order energy-conserving Hamiltonian Boundary Value Methods (HBVMs), and also sketch their efficient implementation for the problem at hand; in Section~\ref{num} we report some numerical tests, showing the potentiality of the proposed numerical approximation procedure; at last, a few conclusions are given in Section~\ref{fine}.

\section{Real form and space semi-discretization}\label{rform}

Let us now write equation (\ref{ms}) in real form, by setting $$\psi(x,t) = u(x,t)+\ii v(x,t),$$ with $u,v:\D\rightarrow \RR^n$. One obtains, by recalling the notation (\ref{v2}):
\begin{equation}\label{rms}
\partial_t u = -\beta \partial_{xx} v -\left[\gamma(u^2+v^2)\right]\circ v, \qquad \partial_t v = \beta \partial_{xx} u+\left[\gamma(u^2+v^2)\right]\circ u, \qquad (x,t)\in\D.
\end{equation}
Similarly, setting $\psi^0 = u^0+\ii v^0$, from (\ref{ms1}) one has that (\ref{rms}) is completed with 
\begin{equation}\label{rms1}
u(x,0) = u^0(x), \qquad v(x,0) = v^0(x), \qquad x\in[a,b].
\end{equation}
and periodic boundary conditions. Consequently, the invariants (\ref{Mj})--(\ref{E}) now become, respectively,
\begin{eqnarray}
M_i(t) &=& \int_a^b \left[u_i(x,t)^2+v_i(x,t)^2\right]\dd x, \qquad i=1,\dots,n, \label{rMj}\\
M(t) &=& \int_a^b \left[ u(x,t)^\top u(x,t)+v(x,t)^\top v(x,t)\right]\dd x, \label{rM}\\
K(t) &=&\int_a^b \left[ \partial_xv(x,t)^\top u(x,t)-\partial_xu(x,t)^\top v(x,t)\right]\dd x, \label{rK}\\ \nonumber
E(t) &=& \frac{1}2\int_a^b \big[ \partial_xu(x,t)^\top\beta\partial_xu(x,t) +  \partial_xv(x,t)^\top\beta\partial_xv(x,t) \\
       && -\frac{1}2(u(x,t)^2+v(x,t)^2)^\top\gamma(u(x,t)^2+v(x,t)^2)\big]\dd x. \label{rE}
\end{eqnarray}
In particular, $E(t)$ coincides with the value of the Hamiltonian functional,
\begin{eqnarray}\nonumber
\H[u,v](t) &=& \frac{1}2\int_a^b \left[\partial_x u^\top\beta \partial_xu+\partial_x v^\top\beta \partial_xv-\frac{1}2(u^2+v^2)^\top\gamma(u^2+v^2)\right]\dd x\\ \label{rH}
&\equiv& \int_a^b \hat H(u,v,\partial_xu,\partial_xv)\dd x,
\end{eqnarray}
which confers a Hamiltonian structure to (\ref{rms}). In fact, one has:
\begin{eqnarray}\nonumber
\partial_t u &=& \delta_v \H[u,v] ~\equiv~ \left(\partial_v-\partial_x(\partial_{\partial_x v})\right)\hat H(u,v,\partial_xu,\partial_xv),\\[-1mm]
\label{rmsH}
\\[-1mm] \nonumber
\partial_t v &=& -\delta_u\H[h,v] ~\equiv~ -\left(\partial_u-\partial_x(\partial_{\partial_x u})\right)\hat H(u,v,\partial_xu,\partial_xv),
\end{eqnarray}
with $\delta_u$ and $\delta_v$ the partial functional derivatives of the functional $\H[u,v]$. 

In order to numerically solve (\ref{rms})-(\ref{rms1}), we consider the expansion of the solution along the periodic orthonormal basis (in space)
\begin{eqnarray}\label{fou}
w_0(x) &=&\frac{1}{\sqrt{b-a}},\\ \nonumber
w_{2j-1}(x) &=&\sqrt{\frac{2}{b-a}}\sin\left( 2\pi j\frac{x-a}{b-a}\right),\qquad x\in[a,b],\\ \nonumber
w_{2j}(x) &=&\sqrt{\frac{2}{b-a}}\cos\left( 2\pi j\frac{x-a}{b-a}\right),\qquad j=1,2,\dots,
\end{eqnarray}
such that, for all allowed indexes,
\begin{equation}\label{orto}
\int_a^b w_i(x)w_j(x)\dd x = \delta_{ij}, 
\end{equation}
with $\delta_{ij}$ the Kronecker delta. In so doing, for suitable time dependent vector coefficients, $q_j(t),p_j(t)\in\RR^n$, one has:
\begin{equation}\label{expuv}
u(x,t) = \sum_{j\ge0} w_j(x)q_j(t),\qquad v(x,t) = \sum_{j\ge0} w_j(x)p_j(t).
\end{equation}
Consequently, the periodic boundary conditions are now always satisfied. Moreover, 
introducing the infinite vector and matrices
\begin{equation}\label{wqp}
\bfw(x) = \pmatrix{c} w_0(x)\\ w_1(x)\\ \vdots\endpmatrix, \quad
\bfq(t) = \pmatrix{ccc} q_0(t),& q_1(t),& \dots\endpmatrix, \quad
\bfp(t) = \pmatrix{ccc} p_0(t),& p_1(t),& \dots\endpmatrix, 
\end{equation}
one has that the expansions (\ref{expuv}) can be, more compactly, written as
\begin{equation}\label{vexpuv}
u(x,t) = \bfq(t)\bfw(x), \qquad v(x,t) = \bfp(t)\bfw(x).
\end{equation}

Introducing the infinite matrices
\begin{equation}\label{D}
D = \frac{2\pi}{b-a}\pmatrix{cccc} 0 \\ &1\cdot I_2\\ && 2\cdot I_2\\ &&& \ddots\endpmatrix, 
\end{equation}
with $I_r\in\RR^{r\times r}$ the identity matrix, and
\begin{equation}\label{tD}
\tD = \frac{2\pi}{b-a}\pmatrix{cccc} 0 \\ &1\cdot J_2\\ && 2\cdot J_2\\ &&& \ddots\endpmatrix, \qquad
J_2=\pmatrix{cc} &1\\-1\endpmatrix, 
\end{equation}
such that
\begin{equation}\label{w12}
\tD^\top=-\tD, \qquad \tD^\top\tD = D^2, \qquad \bfw'(x) = \tD\bfw(x), \qquad \bfw''(x) = \tD^2\bfw(x) \equiv -D^2 \bfw(x),
\end{equation}
one easily proves the following preliminary result.

\begin{lem}\label{deriv}
With reference to the expansions (\ref{expuv}), one has:
\begin{eqnarray*}
u_t(x,t) = \dot\bfq(t)\bfw(x), && v_t(x,t) = \dot\bfp(t)\bfw(x),\\
u_x(x,t)=\bfq(t)\tD\bfw(x), && v_x(x,t)=\bfp(t)\tD\bfw(x),\\
u_{xx}(x,t)=-\bfq(t)D^2\bfw(x), && v_{xx}(x,t)=-\bfp(t)D^2\bfw(x).
\end{eqnarray*}
\end{lem}
\medskip

The following result then holds true.

\begin{theo}\label{rHform} Problem (\ref{rms})-(\ref{rms1}) can be recast  as
\begin{eqnarray}\label{vrms}
\dot\bfq &=& \beta\bfp D^2-\int_a^b\left[\gamma[ (\bfq\bfw(x))^2+(\bfp\bfw(x))^2]\right]\circ (\bfp\bfw(x))\bfw(x)^\top\dd x, \\ \nonumber
\dot\bfp &=& -\beta\bfq D^2+\int_a^b\left[\gamma[ (\bfq\bfw(x))^2+(\bfp\bfw(x))^2]\right]\circ (\bfq\bfw(x))\bfw(x)^\top\dd x, \qquad t\in[0,T],
\end{eqnarray}
with the initial conditions
\begin{equation}\label{vrms1}
\bfq(0) = \bfq^0 \equiv \int_a^b u^0(x)\bfw(x)^\top\dd x, \qquad \bfp(0) = \bfp^0 \equiv \int_a^b v^0(x)\bfw(x)^\top\dd x.
\end{equation}
\end{theo}
\proof
The statement easily follows from Lemma~\ref{deriv}, by right-multiplication by $\bfw(x)^\top$ and then integrating in space, by considering that from (\ref{orto}) one has
\begin{equation}\label{I}
\int_a^b\bfw(x)\bfw(x)^\top\dd x =I.
\end{equation}
with $I$ the identity operator.\,\QED
\bigskip

The following statement concerns the invariants (\ref{rMj})--(\ref{rE}).

\begin{theo}\label{rMj2E}
Setting 
\begin{equation}\label{qjpj}
q_j=\pmatrix{c}q_{1j}\\ \vdots\\ q_{nj}\endpmatrix,\qquad p_j=\pmatrix{c}p_{1j}\\ \vdots\\ p_{nj}\endpmatrix,
\end{equation}
the invariants (\ref{rMj})--(\ref{rE}) can be written as:
\begin{eqnarray}
M_i(t) &=&  \sum_{j\ge0}(q_{ij}^2+p_{ij}^2),  \qquad i=1,\dots,n,\label{rMj1}\\ 
M(t)    &=&  \sum_{j\ge0} \left( q_j^\top q_j+p_j^\top p_j\right) ~\equiv~\sum_{i=1}^n M_i(t),\label{rM1}\\ 
K(t)    &=&  2\sum_{j\ge1} d_{2j}\left( p_{2j}^\top q_{2j-1}-q_{2j}^\top p_{2j-1}\right),\label{rK1}\\
E(t) &=& \frac{1}2\sum_{j\ge0} d_j^2(q_j^\top\beta q_j+p_j^\top\beta p_j) \nonumber\\
&& -\frac{1}4\int_a^b \left[(\bfq(t)\bfw(x))^2+(\bfp(t)\bfw(x))^2\right]^\top\gamma\left[(\bfq(t)\bfw(x))^2+(\bfp(t)\bfw(x))^2\right]\dd x,\label{rE1}
\end{eqnarray}
with 
\begin{equation}\label{dj}
d_j = \frac{2\pi}{b-a}\lceil\frac{j}2\rceil,
\end{equation}
the $j$-th diagonal entry of matrix $D$ in (\ref{D}), $j=0,1,2,\dots$.\footnote{As is usual, $\lceil\cdot\rceil$ denotes the ceiling function.}
\end{theo}
\proof The proof of (\ref{rMj1})  follows from (\ref{vexpuv}) and (\ref{orto}). In fact, denoting $e_i\in\RR^n$  the $i$-th unit vector, one has:
\begin{eqnarray*}
M_i(t) &=& \int_a^b \left( e_i^\top\bfq\bfw\right)^\top\left( e_i^\top\bfq\bfw\right)+
\left( e_i^\top\bfp\bfw\right)^\top\left( e_i^\top\bfp\bfw\right)\dd x \\
&=&\int_a^b  \bfw^\top\left(\bfq^\top e_ie_i^\top\bfq+ \bfp^\top e_i e_i^\top\bfp\right)\bfw\dd x\\
&=&\sum_{k,j\ge0} \left(q_{ik}q_{ij}+p_{ik}p_{ij}\right)\underbrace{\int_a^b w_kw_j\dd x}_{\delta_{kj}}
~=~ \sum_{j\ge0}(q_{ij}^2+p_{ij}^2).
\end{eqnarray*}
Similarly, concerning (\ref{rM1}) one has:
$$
M(t) = \int_a^b \bfw^\top\left(\bfq^\top\bfq+\bfp^\top\bfp\right)\bfw\dd x =\sum_{i,j\ge0} \left(q_i^\top q_j+p_i^\top p_j\right)\underbrace{\int_a^b w_iw_j\dd x}_{=\delta_{ij}} = \sum_{j\ge0} \left(q_j^\top q_j+p_j^\top p_j\right).
$$
The proof of (\ref{rK1}) follows from (\ref{rK}), (\ref{tD}), Lemma~\ref{deriv}, (\ref{orto}), and (\ref{dj}), since
\begin{eqnarray*}
K(t)&=& \int_a^b \bfw(x)^\top\tD\left( \bfq^\top\bfp-\bfp^\top\bfq\right)\bfw(x)\dd x \\
&=&\frac{2\pi}{b-a} \sum_{j\ge1} j\left(p_{2j}^\top q_{2j-1} -p_{2j-1}^\top q_{2j} - q_{2j}^\top p_{2j-1}+q_{2j-1}^\top p_{2j} \right)\\
&=&2 \sum_{j\ge1} d_{2j}\left(p_{2j}^\top q_{2j-1} -q_{2j}^\top p_{2j-1} \right).
\end{eqnarray*}
At last, for (\ref{rE1}), from (\ref{orto}) one has
\begin{eqnarray*}
\int_a^b\partial_x u^\top \beta\partial_x u\dd x &=&-\int_a^b u^\top \beta\partial_{xx} u\dd x ~=~ \int_a^b (\bfq\bfw)^\top\beta(\bfq D^2\bfw)\dd x\\
&=& \int_a^b \bfw^\top\bfq^\top\beta\bfq D^2\bfw\dd x~=~ \sum_{i,j\ge0}  q_i^\top\beta q_jd_j^2\underbrace{\int_a^b w_i(x)w_j(x)\dd x}_{=\delta_{ij}} \\[-2mm]
&=&\sum_{j\ge0} d_j^2 q_j^\top\beta q_j,
\end{eqnarray*} 
and, similarly, one obtains
$$\int_a^b\partial_x v^\top \beta\partial_x v\dd x ~=~ \sum_{j\ge0} d_j^2 p_j^\top\beta p_j.$$
\QED
\bigskip

As a straightforward consequence, the following result holds true.

\begin{cor}\label{rHam}
The infinite-dimensional ODE problem (\ref{vrms}) can be cast in Hamiltonian form as
$$\dot\bfq = \partial_\bfp H(\bfq,\bfp), \qquad \dot\bfp = -\partial_\bfq H(\bfq,\bfp),$$
w.r.t. the Hamiltonian (see (\ref{dj}))
\begin{eqnarray}\nonumber
H(\bfq,\bfp) &=&  \frac{1}2 \sum_{j\ge0} d_j^2\left(q_j^\top\beta q_j+p_j^\top\beta p_j\right) \\ \label{rH1}
&& -\frac{1}4\int_a^b \left[(\bfq(t)\bfw(x))^2+(\bfp(t)\bfw(x))^2\right]^\top\gamma\left[(\bfq(t)\bfw(x))^2+(\bfp(t)\bfw(x))^2\right]\dd x.
\end{eqnarray} 
This latter is equivalent to the Hamiltonian functional (\ref{rH}).
\end{cor}

\subsection{Space semi-discretization}\label{ssd}

As is clear, in order for numerically solving (\ref{vrms})-(\ref{vrms1}), the infinite expansions (\ref{expuv}) need to be truncated to finite sums, at  a convenient index $2N$. In so doing, the infinite vectors and matrices (\ref{wqp})--(\ref{tD}) now respectively become\,\footnote{In order not to make the notation cumbersome, we continue to use the same identifiers in (\ref{wqp})--(\ref{tD}), even though, hereafter, they will denote the finite counterparts.}
\begin{eqnarray}\label{wqpN}
\bfw(x) &=& \pmatrix{c} w_0(x)\\ \vdots \\ w_{2N}(x)\endpmatrix \in\RR^{2N+1},\\ \nonumber
\bfq(t) &=& \pmatrix{ccc} q_0(t),& \dots, &q_{2N}(t)\endpmatrix, ~~
\bfp(t) ~=~ \pmatrix{ccc} p_0(t),& \dots, &p_{2N}(t)\endpmatrix~\in~\RR^{n\times(2N+1)}, 
\end{eqnarray}
and
\begin{equation}\label{DtDN}
D = \frac{2\pi}{b-a}\pmatrix{cccc} 0 \\ &1\cdot I_2\\ &&\ddots\\ &&&  N\cdot I_2\endpmatrix, ~
\tD = \frac{2\pi}{b-a}\pmatrix{cccc} 0 \\ &1\cdot J_2\\ && \ddots\\ &&&N\cdot J_2\endpmatrix \in\RR^{(2N+1)\times (2N+1)}. 
\end{equation}
In so doing, the vector form of the expansions (\ref{vexpuv}), the properties (\ref{w12}), and the result of Lemma~\ref{deriv} continue formally to hold. 
However, in this case, the approximations (\ref{vexpuv}) will no more satisfy, in general, the equations (\ref{rms}). Nevertheless, in the spirit of Galerkin methods, by requiring that each component of the residual be orthogonal to the functional space (see (\ref{fou}))
$${\cal V}_N = \mathrm{span}\left\{w_0(x),w_1(x),\dots,w_{2N}(x)\right\},$$
to which the entries of the approximations (\ref{vexpuv})  belong for all $t$, one obtains that the result of Theorem~\ref{rHform}  continues formally to hold. In particular, the dimension of the (finite-dimensional) Hamiltonian ODE problem (\ref{vrms}) is now $2n(2N+1)$. Moreover, Theorem~\ref{rMj2E} now reads as follows.

\begin{theo}\label{rMj2EN}
Setting $q_j$ and $p_j$ as in (\ref{qjpj}) the invariants (\ref{rMj})--(\ref{rE}) can be written as:
\begin{eqnarray}
M_i(t) &=&  \sum_{j=0}^{2N}(q_{ij}^2+p_{ij}^2),  \qquad i=1,\dots,n,\label{rMj1N}\\ 
M(t)    &=&  \sum_{j=0}^{2N} \left( q_j^\top q_j+p_j^\top p_j\right) ~\equiv~\sum_{i=1}^n M_i(t),\label{rM1N}\\ 
K(t)    &=&  2\sum_{j=1}^{N} d_{2j}\left( p_{2j}^\top q_{2j-1}-q_{2j}^\top p_{2j-1}\right),\label{rK1N}\\
E(t) &=& \frac{1}2\sum_{j=0}^{2N} d_j^2(q_j^\top\beta q_j+p_j^\top\beta p_j) \nonumber\\
&& -\frac{1}4\int_a^b \left[(\bfq(t)\bfw(x))^2+(\bfp(t)\bfw(x))^2\right]^\top\gamma\left[(\bfq(t)\bfw(x))^2+(\bfp(t)\bfw(x))^2\right]\dd x,\label{rE1N}
\end{eqnarray}
with $d_j$ defined as in (\ref{dj}).
\end{theo}

Similarly, Corollary~\ref{rHam} is modified as reported below.

\begin{cor}\label{rHamN}
The finite-dimensional problem (\ref{vrms}) can be cast in Hamiltonian form as
$$\dot\bfq = \partial_\bfp H(\bfq,\bfp), \qquad \dot\bfp = -\partial_\bfq H(\bfq,\bfp),$$
w.r.t. the Hamiltonian (see (\ref{dj}))
\begin{eqnarray}\nonumber
H(\bfq,\bfp) &=&  \frac{1}2 \sum_{j=0}^{2N} d_j^2\left(q_j^\top\beta q_j+p_j^\top\beta p_j\right) \\ \label{rH1N}
&& -\frac{1}4\int_a^b \left[(\bfq\bfw(x))^2+(\bfp\bfw(x))^2\right]^\top\gamma\left[(\bfq\bfw(x))^2+(\bfp\bfw(x))^2\right]\dd x.
\end{eqnarray} 
\end{cor}

It is well-known (see, e.g., \cite{B2001}) that, under regularity assumptions, the truncated versions of the solution, of the invariants, and of the Hamiltonian converge  exponentially to the exact counterparts as $N\rightarrow\infty$ (this fact is often referred to as {\em spectral accuracy}). In so doing, one obtains a Fourier-Galerkin space semi-discretization of the original problem (\ref{rms}).

\section{Hamiltonian Boundary Value Methods}\label{hbvms}

For conveniently introducing the time integration of system (\ref{vrms})-(\ref{vrms1}), let us define matrix
\begin{equation}\label{y}
\bfy \equiv \pmatrix{ccc} y_0,&\dots,&y_{2N}\endpmatrix\in\RR^{2n\times(2N+1)},
\end{equation}
with (see (\ref{qjpj}))
\begin{equation}\label{yj}
y_j = \pmatrix{ccccccc}q_{1j}, &p_{1j}, &q_{2j}, &p_{2j}, &\dots, &q_{nj}, &p_{nj}\endpmatrix^\top\in\RR^{2n}, \quad j=0,\dots,2N,
\end{equation}
and the block diagonal matrices (see (\ref{tD})) 
\begin{equation}\label{b2g2}
\beta_2 = \beta\otimes I_2, \qquad \gamma_2 = \gamma\otimes I_2, \qquad J = I_n\otimes J_2, \qquad Q = I_n\otimes 
\pmatrix{cc} 1&1\\1&1\endpmatrix.
\end{equation}
In so doing, the system of ODEs (\ref{vrms}) can be rewritten as (see (\ref{DtDN}))
\begin{equation}\label{rmsy}
\dot\bfy = J\left[\beta_2 \bfy D^2 -\int_a^b \left[\gamma_2Q\left(\bfy\bfw(x)\right)^2\right]\circ\left(\bfy\bfw(x)\right)\bfw(x)^\top\dd x\right] \equiv J\nabla H(\bfy), \quad t\in[0,T],
\end{equation}
which is Hamiltonian w.r.t. the Hamiltonian (see (\ref{rH1N}))
\begin{equation}\label{H1Ny}
H(\bfy) = \frac{1}2\sum_{j=0}^{2N} d_j^2\left(y_j^\top\beta_2 y_j\right)-\frac{1}8\int_a^b \left[ Q\left(\bfy\bfw(x)\right)^2\right]^\top\gamma_2\left[ Q\left(\bfy\bfw(x)\right)^2\right]\dd x.
\end{equation}
Similarly, the initial condition (\ref{vrms1}) becomes
\begin{equation}\label{rms1y}
\bfy(0) = \bfy^0\equiv \pmatrix{ccccccc}q_{1j}^0, &p_{1j}^0, &q_{2j}^0, &p_{2j}^0, &\dots, &q_{nj}^0, &p_{nj}^0\endpmatrix^\top,
\end{equation}
where the same notation (\ref{qjpj}) has been used for the columns of $\bfq^0$ and $\bfp^0$ in (\ref{vrms1}). 
At last, the invariants (\ref{rMj1N})--(\ref{rE1N}) respectively read, by setting $y_{ij}$ the $i$-th entry of $y_j$, and with reference to (\ref{H1Ny}),
\begin{eqnarray}
M_i(t) &=&   \sum_{j=0}^{2N} (y_{2i-1,j}^2+y_{2i,j}^2),  \qquad i=1,\dots,n, \label{Mjy}\\
M(t)   &=&    \sum_{j=0}^{2N} y_j^\top y_j \,\equiv\, \sum_{i=1}^n M_i(t),  \label{My}\\
K(t)   &=&      2\sum_{j=1}^N d_{2j}\,y_{2j-1}^\top J y_{2j},\label{Ky}\\
E(t)   &\equiv& H(\bfy(t)). \label{Ey}
\end{eqnarray}

With these premises, we shall consider the time integration of the Hamiltonian ODE problem (\ref{rmsy})-(\ref{rms1y}) by using HBVM$(k,s)$ methods, where HBVM is the acronym of {\em Hamiltonian Boundary Value Methods} and the couple of parameters $(k,s)$ characterizes the specific method. Such methods, developed in \cite{LIMbook2016,BIT2009_1,BIT2010,BIT2012,BIT2015} (see also the recent review paper \cite{BI2018}) are a class of energy-conserving Runge-Kutta (RK, hereafter) methods for Hamiltonian problems. HBVMs have been generalized along several directions, including Hamiltonian BVPs \cite{ABI2015}, Poisson problems \cite{BCMR2012,BMR2019}, problem with multiple invariants \cite{BI2012,BS2014}, constrained Hamiltonian problems \cite{BGIW2018}, fractional differential equations \cite{ABI2019_1}, highly oscillatory problems \cite{BMR2018}. They have been also used as spectral methods in time \cite{ABI2019_1,ABI2018} and for the time integration of Hamiltonian PDEs \cite{BBFCI2018,BFCI2015,BFCI2019,BGS2019,BGZ2019,LIMbook2016,BZL2018}. The efficient implementation of HBVMs has been studied in \cite{BFCI2014,BIT2009,BIT2011} (see also \cite[Chapter\,4]{LIMbook2016}). Such methods also admit a continuos-stage RK form \cite{ABI2019_1}. In more detail, a HBVM$(k,s)$ method is the $k$-stage RK method, whose Butcher tableau is given by
\begin{equation}\label{hbvmks}
\begin{array}{c|c} c & \I_s \P_s^\top\Omega \\ \hline\\[-3mm] & b^\top\end{array} ~,
\end{equation}
with
\begin{equation}\label{bc}
c = \pmatrix{ccc} c_1,&\dots,&c_k\endpmatrix^\top, \qquad b = \pmatrix{ccc} b_1,&\dots,&b_k\endpmatrix^\top,
\end{equation}
the abscissae and weights of the Gauss-Legendre formula of order $2k$, and
\begin{eqnarray}\nonumber
\P_s &=& \pmatrix{ccc} P_0(c_1) & \dots & P_{s-1}(c_1)\\
\vdots & &\vdots\\
P_0(c_k) & \dots & P_{s-1}(c_k)
\endpmatrix,\qquad \Omega = \pmatrix{ccc} b_1\\ &\ddots\\ &&b_k\endpmatrix,\\[2mm]
 \I_s &=& \pmatrix{ccc} \int_0^{c_1}P_0(x)\dd x & \dots & \int_0^{c_1} P_{s-1}(x)\dd x\\
\vdots & &\vdots\\
\int_0^{c_k}P_0(x)\dd x & \dots & \int_0^{c_k} P_{s-1}(x)\dd x
\endpmatrix.\label{IPO}
\end{eqnarray}
The following result summarizes some of the properties of HBVMs \cite{LIMbook2016,BI2018,BIT2012}.

\begin{theo}\label{hbvmth}
For all $k\ge s$, a HBVM$(k,s)$ method:
\begin{itemize}
\item is symmetric and has order $2s$;
\item coincides with the (symplectic) $s$-stage Gauss collocation method, when $k=s$;
\item is energy conserving for all polynomial Hamiltonians of degree not larger than $2k/s$.
\end{itemize}
\end{theo}
Since the Hamiltonian (\ref{H1Ny}) is a polynomial of degree 4 in $\bfy$, the following result easily follows.

\begin{cor}\label{ep}
For  all $s\ge1$, HBVM$(2s,s)$ methods are energy-conserving and of order $2s$, when applied for solving problem (\ref{rmsy})-(\ref{rms1y}).\footnote{In particular, when $s=1$ one retrieves the AVF method in \cite{avf} applied for solving (\ref{rmsy})-(\ref{rms1y}).}
\end{cor}

\subsection{Solving the discrete problem}\label{sdp}
Let us now sketch the efficient solution of the discrete problem generated by a HBVM$(k,s)$ method, with $k\ge s$, applied for solving (\ref{rmsy})-(\ref{rms1y}). The stage-equation of the first step of integration, when using a stepsize $h$, reads
\begin{equation}\label{Y}
Y = e\otimes\bfy^0 + h\I_s\P_s^\top\Omega\otimes J\, \nabla H(Y),
\end{equation}
where $e=\pmatrix{ccc} 1,& \dots, &1\endpmatrix^\top\in\RR^{k}$, and
$$Y=\pmatrix{c} Y_1\\ \vdots\\ Y_{k}\endpmatrix,\, \nabla H(Y)=\pmatrix{c}\nabla H(Y_1)\\ \vdots\\ \nabla H(Y_{k})\endpmatrix
\in~\RR^{2kn\times (2N+1)}, $$
with$$Y_j,\,\nabla H(Y_j)\in\RR^{2n\times (2N+1)}, \qquad j=1,\dots,k.$$
In order for reducing the size of the discrete problem, according to \cite{BIT2011}, let us set
$$
\Gamma = \P_s^\top\Omega\otimes J\, \nabla H(Y) ~\in~\RR^{2sn\times(2N+1)}.
$$
so that (\ref{Y}) reads
$$Y = e\otimes\bfy^0 + h\I_s\otimes I_{2n}\,\Gamma.$$
From the last two equations, one obtains, at last,
\begin{equation}\label{Gam}
\Gamma \equiv \pmatrix{c} \Gamma_0\\ \vdots\\ \Gamma_{s-1}\endpmatrix\,=\, \P_s^\top\Omega\otimes J\, \nabla H\left(e\otimes\bfy^0 + h\I_s\otimes I_{2n}\,\Gamma\right),\end{equation}
with
\begin{equation}\label{Gamj}
\Gamma_j\in\RR^{2n\times(2N+1)}, \qquad j=0,\dots,s-1.
\end{equation}
According to \cite{BIT2011}, the new approximation is then given by 
\begin{equation}\label{y1}
\bfy^1 = \bfy^0 + h\,\Gamma_0.
\end{equation}
Consequently, the discrete problem (\ref{Gam}) has now dimension $2sn\times (2N+1)$, and we need to solve the nonlinear matrix equation
\begin{equation}\label{mateq} 
\Gamma-\P_s^\top\Omega\otimes J\, \nabla H\left(e\otimes\bfy^0 + h\I_s\otimes I_{2n}\,\Gamma\right) = O.
\end{equation}
We observe that the above problem always admits a solution, provided that the stepsize $h$ is small enough. In order to discuss this issue without too many complicated details, we shall consider the simpler case where, in (\ref{b2g2})-(\ref{rmsy}), $\gamma=0$. In fact, in such a case, one derives necessary conditions for the existence of the solution of the general problem, even though $\|\gamma_2\|\ll\|D^2\|$, when $N\gg1$. The following result then holds true.

\begin{theo}\label{hsmall} Assume that in (\ref{b2g2})-(\ref{rmsy}) $\gamma=0$. Then, for all stepsizes $h>0$ satisfying 
\begin{equation}\label{hsmall1} h<\|\beta\|^{-1}\left(\frac{b-a}{2N}\right)^2\end{equation}
the matrix problem (\ref{mateq}) for the unknown $\Gamma$ is solvable and the fixed-point iteration
\begin{equation}\label{fixit}
\Gamma^{\ell+1} = \P_s^\top\Omega\otimes J\, \nabla H\left(e\otimes\bfy^0 + h\I_s\otimes I_{2n}\,\Gamma^\ell\right), \qquad \ell=0,1,\dots,
\end{equation}
converges to its solution.
\end{theo}
\proof
In fact, by taking into account that
\begin{equation}\label{Xs}
X_s \,\equiv\, \P_s^\top\Omega\I_s = \pmatrix{cccc} \xi_0 &-\xi_1\\ \xi_1 & 0 &\ddots\\ &\ddots &\ddots &-\xi_{s-1}\\ &&\xi_{s-1} &0\endpmatrix,
\quad \xi_i = \left(2\sqrt{|4i^2-1|}\right)^{-\frac{1}2}, \quad i=0,\dots,s-1.
\end{equation} 
and ~$\P_s^\top\Omega e=(\,1,\,0,\,\dots,\,0\,)^\top\equiv e_1^\top\in\RR^s$,  the fixed-point iteration (\ref{fixit}) in this case reads
$$\Gamma^{\ell+1} = \left[e_1\otimes (\beta\otimes J_2\, \bfy^0) + h(X_s\otimes\beta\otimes J_2)\,\Gamma^\ell\right]D^2,\qquad \ell=0,1,\dots.$$
Consequently, the second member of the equation turns out to be a contraction, provided that
$$h\|X_s\|\cdot\|\beta\|\cdot\|D^2\|<1.$$ Eventually, the statement follows by considering that (see (\ref{DtDN}) and (\ref{Xs}))
$$\|X_s\|<1,\qquad \mbox{and}\qquad \|D^2\|=\left(\frac{2N}{b-a}\right)^2.\,\QED$$

\begin{rem} As is clear, the bound (\ref{hsmall1}) prevents the usage of relatively large  time-steps when a spectrally accurate space discretization is considered (so that $N\gg1$), and the (straightforward) fixed-point iteration (\ref{fixit}) implements the solution of (\ref{mateq}).\end{rem}

As stated in the previous remark, if the usage of large stepsizes is sought, as in the case of spectral methods in time, we need to resort to a Newton-type iteration for solving (\ref{mateq}). For this purpose, let us now recall the ``$\vec$'' function
$$\vec\left[\pmatrix{ccc} a_{11} & \dots &a_{1n}\\\vdots & &\vdots\\ a_{m1}& \dots&a_{mn}\endpmatrix\right] = \pmatrix{c}
a_{11}\\ \vdots \\ a_{m1}\\ \vdots \\ a_{1n}\\ \vdots\\ a_{mn}\endpmatrix\in\RR^{m\cdot n},$$
and the property, with $A,X,B$ suitable given matrices, $$\vec(AXB) = \left(B^\top\otimes A\right)\vec(X).$$ 
In order to solve (\ref{mateq}), it is then convenient to define the equivalent vector of the unknown
\begin{equation}\label{vgam}
\bfg = \vec\left(\Gamma^\top\right) \in\RR^{s\,2n\,(2N+1)}.
\end{equation}
In terms of such unknown vector, the matrix problem (\ref{mateq}) becomes
\begin{equation}\label{vmateq}
\bff(\bfg) \,\equiv\, \bfg -  P_s^\top\Omega\otimes J \otimes I_{2N+1}\,\vec\left[\nabla H\left(e\otimes\bfy^0+h\I_s\otimes I_{2n}\,\Gamma\right)^\top\right] = \bfzero.
\end{equation}
Let us hereafter set, for the sake of brevity, 
\begin{equation}\label{I}
I\equiv I_{s\,2n\,(2N+1)}.
\end{equation} 
Moreover, we observe that, for enough regular and bounded solutions of (\ref{ms})-(\ref{ms1}) (which we have assumed), one has that the linear term at the right-hand side in (\ref{rmsy}) dominates, when $N\gg1$ (as is the case, when a spectrally accurate space semi-discretization is sought).
Consequently, the simplified Newton iteration for solving (\ref{vmateq}), with the Hessian approximated by the linear part alone reads, by taking into account (\ref{b2g2}) and (\ref{Xs}),
\begin{eqnarray}\nonumber
\mathrm{solve:} && \left[ I-h X_s\otimes \beta\otimes J_2\otimes D^2\right]\bfdelta^\ell = -\bff(\bfg^\ell),\\
\mathrm{set:}&& \bfg^{\ell+1} = \bfg^\ell +\bfdelta^\ell, \qquad \ell=0,1,\dots.\label{simpNewt}
\end{eqnarray}
This straightforward iteration has the obvious advantage of having a constant coefficient matrix for the linear systems to be solved at all time steps. Nevertheless, its dimension may be quite large, i.e., $s\, 2n\, (2N+1)$. Nevertheless, a corresponding {\em blended iteration}, having a much more favourable complexity, can be considered. This iteration, at first devised in \cite{BM2002,BM2009} for block implicit methods, has been implemented in the computational codes {\tt BiM} \cite{BM2004}, for stiff ODE-IVPs,  and {\tt BiMD} \cite{BMM2006}, also solving DAEs. Later on, it has been considerd for HBVMs \cite{BFCI2014,LIMbook2016,BIT2011} and, more recently, for RKN-type methods \cite{WMF2017}. It is worth mentioning that, in the case of HBVMs, it has allowed their usage as spectral methods in time \cite{ABI2019,BIMR2018,BMR2018}, because the use of relatively large stepsizes has been made possible.
In order to provide some detail on this iteration, let us consider the solution of the linear system
\begin{equation}\label{uno}
\left[ I-h X_s\otimes \beta\otimes J_2\otimes D^2\right]\bfdelta = -\bff(\bfg) \,\equiv\, \bfeta
\end{equation}
associated with (\ref{simpNewt}), where we have omitted the superscript $\ell$, for the sake of brevity. Since matrix $X_s$ is known to be nonsingular, left multiplication of both members of (\ref{uno}) by $$\rho_sX_s^{-1}\otimes I_{2n(2N+1)},$$ with $\rho_s$ a positive parameter to be chosen later,  allows to obtain the equivalent linear system
\begin{equation}\label{due}
\rho_s\left[ X_s^{-1}\otimes I_{2n(2N+1)}-h I_s\otimes \beta\otimes J_2\otimes D^2\right]\bfdelta = \left[\rho_sX_s^{-1}\otimes I_{2n(sN+1)}\right] \bfeta \,\equiv\, \bfeta_1.
\end{equation}
Next, by defining the {\em weighting function}
\begin{equation}\label{teta}
\Theta = \left( I-h\rho_s I_s\otimes \beta\otimes J_2\otimes D^2\right)^{-1},
\end{equation}
we obtain a further equivalent system as the {\em blending} of (\ref{uno}) and (\ref{due}), with weights $\Theta$ and $I-\Theta$, respectively. Consequently,  by calling $M$ the resulting coefficient matrix, one  obtains the linear system
$$M\bfdelta = \bfeta_1+\Theta(\bfeta-\bfeta_1).$$
This latter system is then solved by considering the splitting
$$M = N-(N-M), \qquad N = \left( I-h\rho_s I_s\otimes \beta\otimes J_2\otimes D^2\right)\,\equiv\,\Theta^{-1}.$$
The corresponding {\em blended iteration} is then given by
$$\bfdelta^0=\bfzero,\qquad \bfdelta^\ell = \Theta\left[ (N-M)\bfdelta^{\ell-1}+ \bfeta_1+\Theta(\bfeta-\bfeta_1)\right], \quad \ell=1,2,\dots.$$
The positive parameter $\rho_s$ is then chosen in order to optimize the convergence properties of the iteration, according to a linear analysis of convergence \cite{BM2002,BM2009}, and turns out to be given by
\begin{equation}\label{ros}
\rho_s = \min_{\lambda\in\sigma(X_s)}|\lambda|,
\end{equation}
where $\sigma(X_s)$ denotes the set of the eigenvalues of $X_s$. Finally, by using a splitting-Newton iteration \cite{BIT2011,BM2004,BMM2006}, one eventually obtains that (\ref{simpNewt}) is replaced by:
\begin{eqnarray}\nonumber
\mathrm{set:} && \bfeta^\ell = -\bff(\bfg^\ell),\\ \nonumber
\mathrm{set:} && \bfeta_1^\ell = \left[\rho_sX_s^{-1}\otimes I_{2n(2N+1)}\right]\bfeta^\ell,\\ \nonumber
\mathrm{compute:} && \bfdelta^\ell = \Theta\left[ \bfeta_1^\ell+\Theta(\bfeta^\ell-\bfeta_1^\ell)\right],\\
\mathrm{set:}&& \bfg^{\ell+1} = \bfg^\ell +\bfdelta^\ell, \qquad\qquad \ell=0,1,\dots.\label{blend}
\end{eqnarray}
Consequently, the complexity of the iteration essentially amounts to the computation of matrix $\Theta$ in (\ref{teta}) which, however, needs to be done only once, since it is the same for all time-steps. Next result, which generalizes that obtained in \cite[Theorem\,5]{BBFCI2018} for the NLSE, shows that such a matrix is block diagonal, with the diagonal blocks having a block-diagonal structure.

\begin{theo}\label{tetath} With reference to (\ref{gbp}), (\ref{DtDN}), and (\ref{ros}), matrix $\Theta$ defined in (\ref{teta}) turns out to be always nonsingular, and is given by
$$\Theta = I_s\otimes \pmatrix{ccc} \theta_1\\ & \ddots \\ &&\theta_n\endpmatrix,$$
with
\begin{equation}\label{tetai}
\theta_i = \pmatrix{cc}
(I_{2N+1}+B_i^2)^{-1} & B_i(I_{2N+1}+B_i^2)^{-1}\\ -B_i(I_{2N+1}+B_i^2)^{-1} &(I_{2N+1}+B_i^2)^{-1}\endpmatrix, \qquad B_i = h\rho_s\beta_iD^2, \qquad i = 1,\dots,n.
\end{equation}
\end{theo}
\proof In fact, from (\ref{I}), (\ref{teta}), and (\ref{tetai}) one obtains that
$$\Theta = I_s\otimes\pmatrix{ccc} \left(I_{2(2N+1)}-J_2\otimes B_1\right)\\
&\ddots\\ && \left(I_{2(2N+1)}-J_2\otimes B_n\right)\endpmatrix^{-1}.$$ 
The statement then follows by considering that
$$\left(I_{2(2N+1)}-J_2\otimes B_i\right)^{-1}\,\equiv\,\pmatrix{cc}
I_{2N+1} &-B_i\\ B_i & I_{2N+1}\endpmatrix^{-1} =
\theta_i,\qquad i=1,\dots,n,$$
and, moreover, such matrices are always well-defined, since (see (\ref{tetai})) $I_{2N+1}+B_i^2$ is a diagonal matrix with positive diagonal entries.\,\QED
\bigskip

For sake of completeness, we also mention that, in order to compute the right-hand side in (\ref{rmsy}) and the Hamiltonian (\ref{H1Ny}), one needs to evaluate the involved integrals in space. Since the integrand is a trigonometric polynomial of degree at most $4N$, one has that they can be exactly computed by using a composite trapezoidal rule at the evenly spaced points
\begin{equation}\label{xi}
x_i = a + i\frac{b-a}m, \qquad i=0,1,\dots,m,
\end{equation}
with $m\ge 4N+1$ \cite{G1997}. Consequently, the choice $m=4N+1$ will be always considered, in the numerical tests. It is worth mentioning that this procedure can be performed with complexity $O(nN\log N)$ by using FFT \cite{BGS2019}, when large values of $N$ are used,  even though we shall not consider it here, for the numerical tests (and, thus, the complexity turns out to be $O(nN^2)$).

\subsection{Spectral HBVMs}\label{spectrem} 
It is worth mentioning that, by choosing $s$ large enough, one can use a HBVM$(k,s)$ method as a spectral method in time. This, in turn,  allows considering relatively large time-steps. The use of HBVMs in this fashion has been considered in \cite{BMR2018} for highly oscillatory problems, and in \cite{BIMR2018} for deriving a spectrally accurate space-time numerical solution of some Hamiltonian PDEs (see also \cite{BIMR2019}). A thorough convergence analysis of the methods, when used as spectral methods in time, has been made in \cite{ABI2019}. Such an analysis allows to state that, when using a finite precision arithmetic with machine epsilon $u$, a spectral accuracy in time can be expected when, with reference to (\ref{Gam}),\footnote{We refer to \cite{ABI2019} for full details.}
$$\|\Gamma_{s-1}\|<tol\cdot \max_{i<s}\|\Gamma_i\|, \qquad \mbox{with}\qquad tol\sim\sqrt{u}.$$ Moreover,
when using a double precision IEEE arithmetic, the value of $k$ can be conveniently chosen by means of the following heuristics \cite{BIMR2018,BMR2018}
$$k = \max\{20,s+2\}.$$
It must be stressed, however, that the usage of a HBVM$(k,s)$ method as a spectral method in time is meaningful only if we can efficiently solve the generated discrete problems, even when a relatively large time-step is considered: this is made possible by the {\em blended iteration} (\ref{blend}).

\section{Numerical tests}\label{num}

We here report a couple of numerical tests, in order to assess the theoretical achievements of the previous sections. All the numerical tests have been performed on an Intel core I7 laptop with 8GB of memory, using Matlab (R2019a). HBVMs for the Manakov problem are implemented by suitably modifying the Matlab function {\tt hbvm} available at \cite{hbvm}. The blended iteration (\ref{blend}), with an approximate Jacobian given by the linear part alone, has been used for solving the discrete problems.

\subsection*{First test problem}  
\begin{figure}[t]
\centerline{\includegraphics[width=5.5cm,height=5cm]{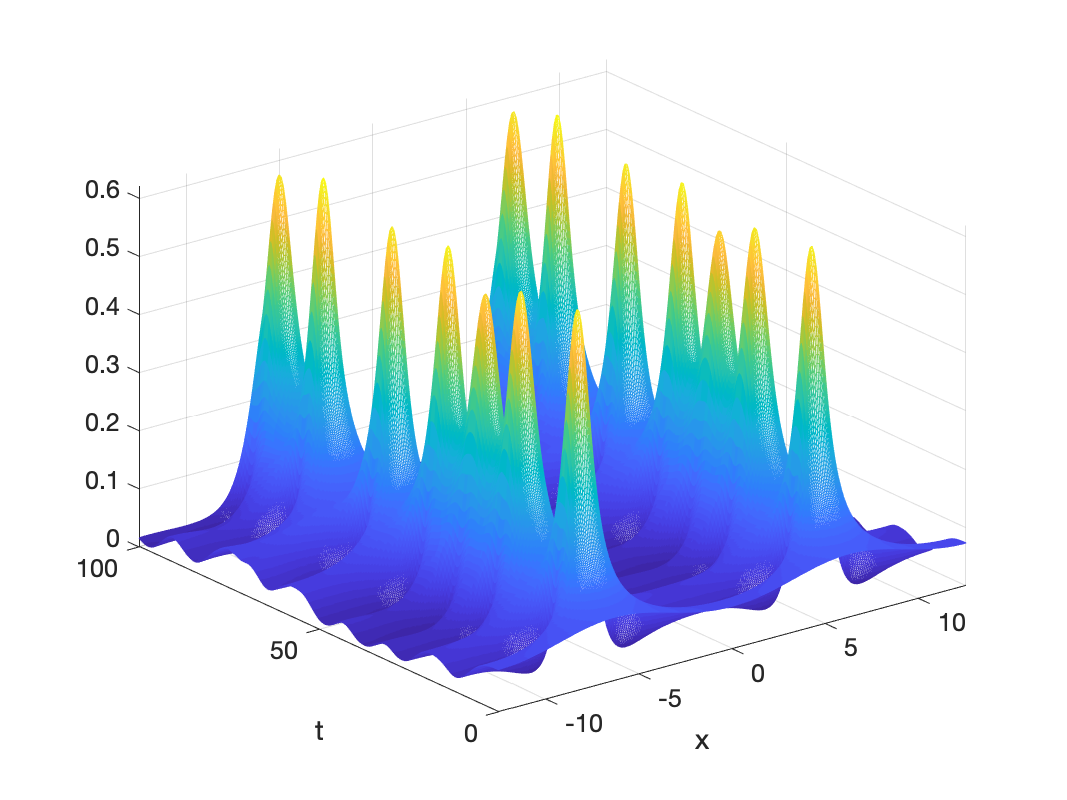}\includegraphics[width=5.5cm,height=5cm]{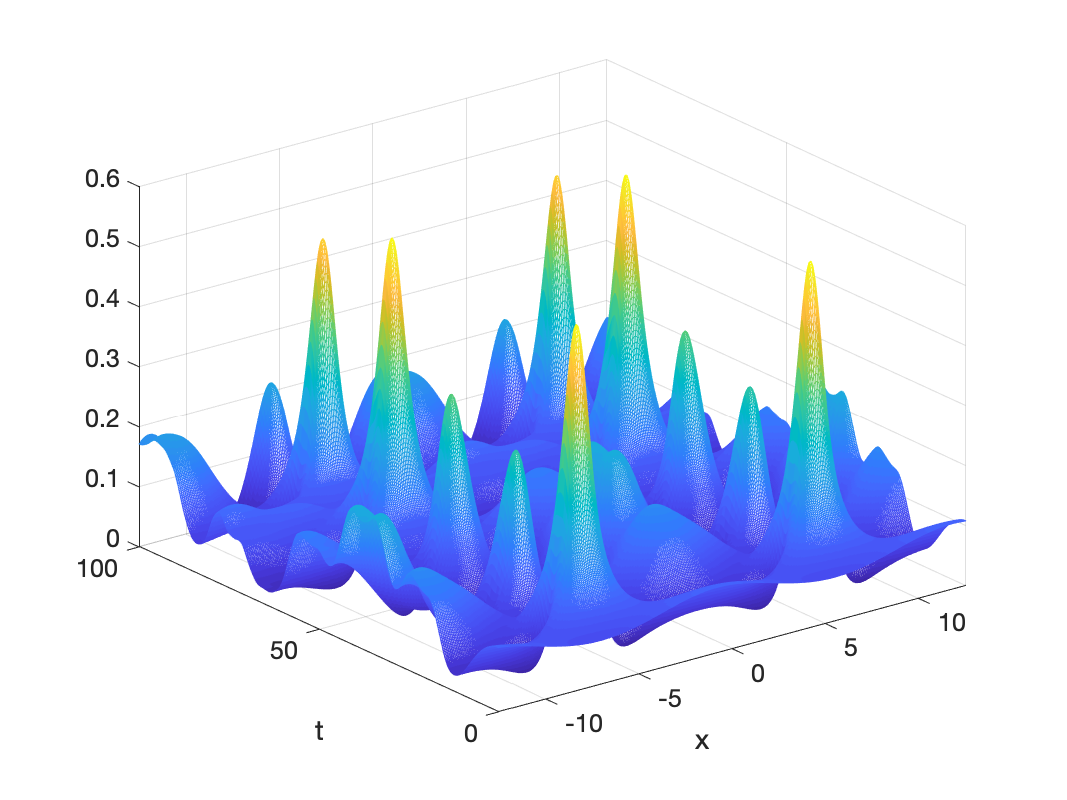}
\includegraphics[width=5.5cm,height=5cm]{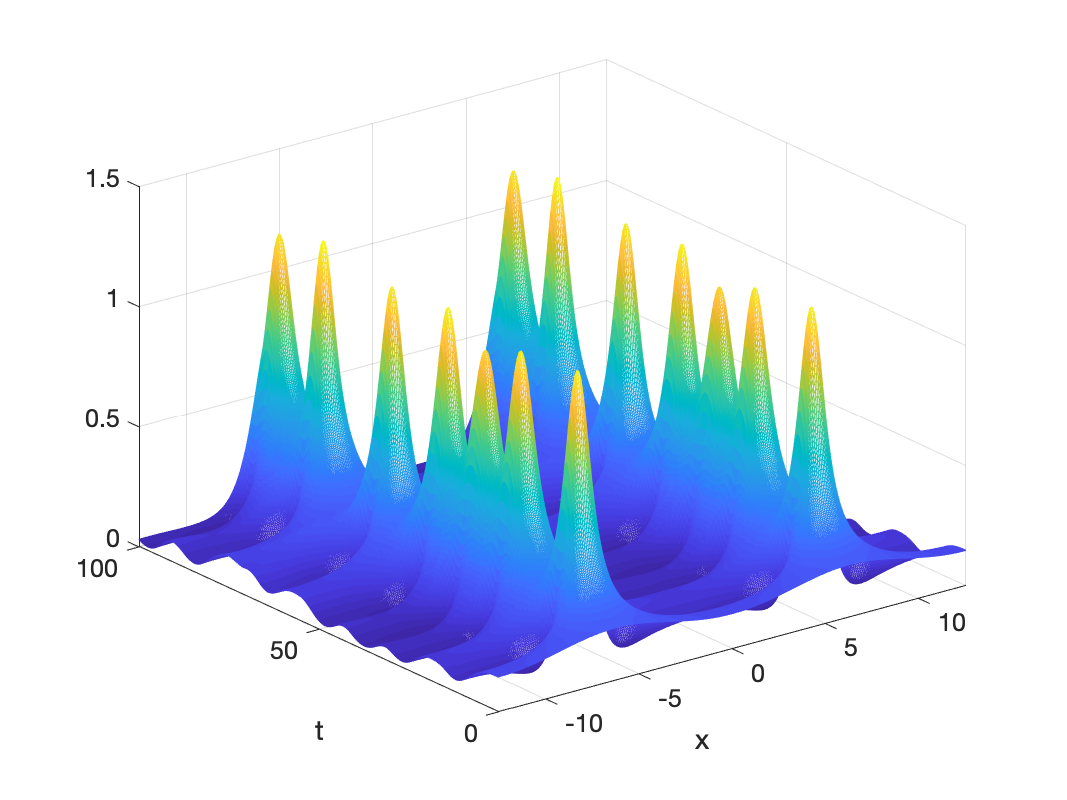}}
\caption{Plots of $|\psi_i(x,t)|^2$, $i=1,2,3$, for problem (\ref{p1})-(\ref{p1par}).}
\label{solu}
\end{figure}

\begin{figure}[t]
\centerline{\includegraphics[width=8.5cm,height=6cm]{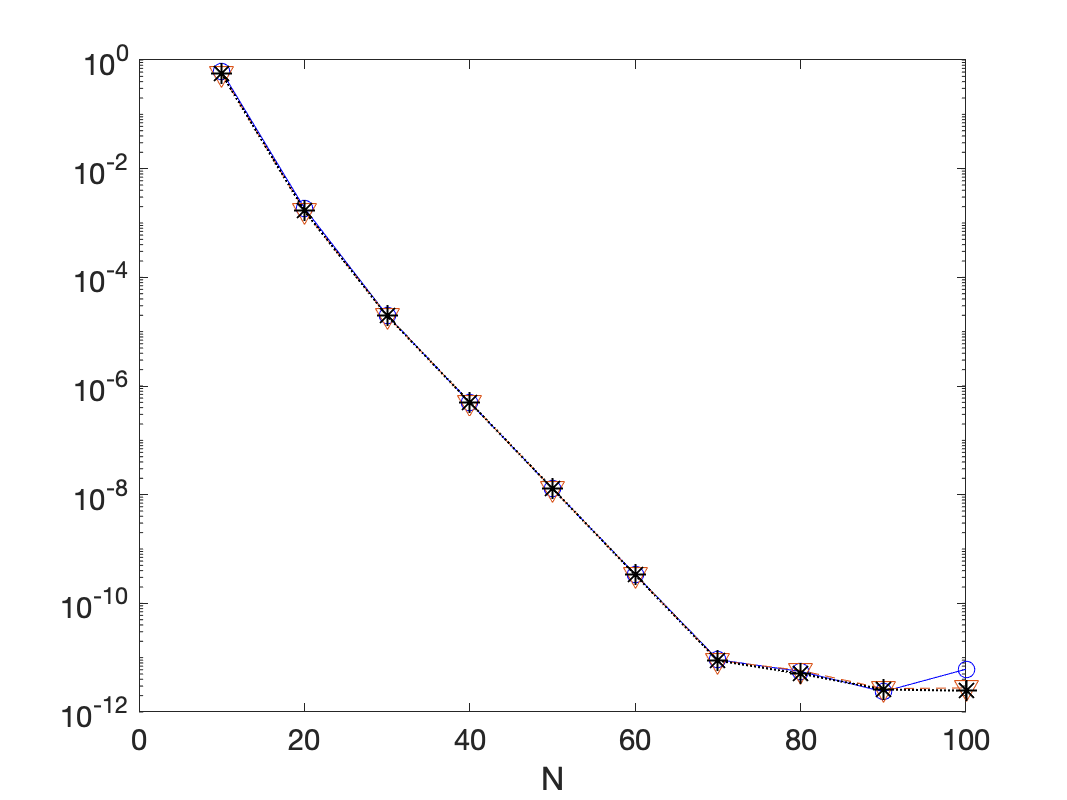}\includegraphics[width=8.5cm,height=6cm]{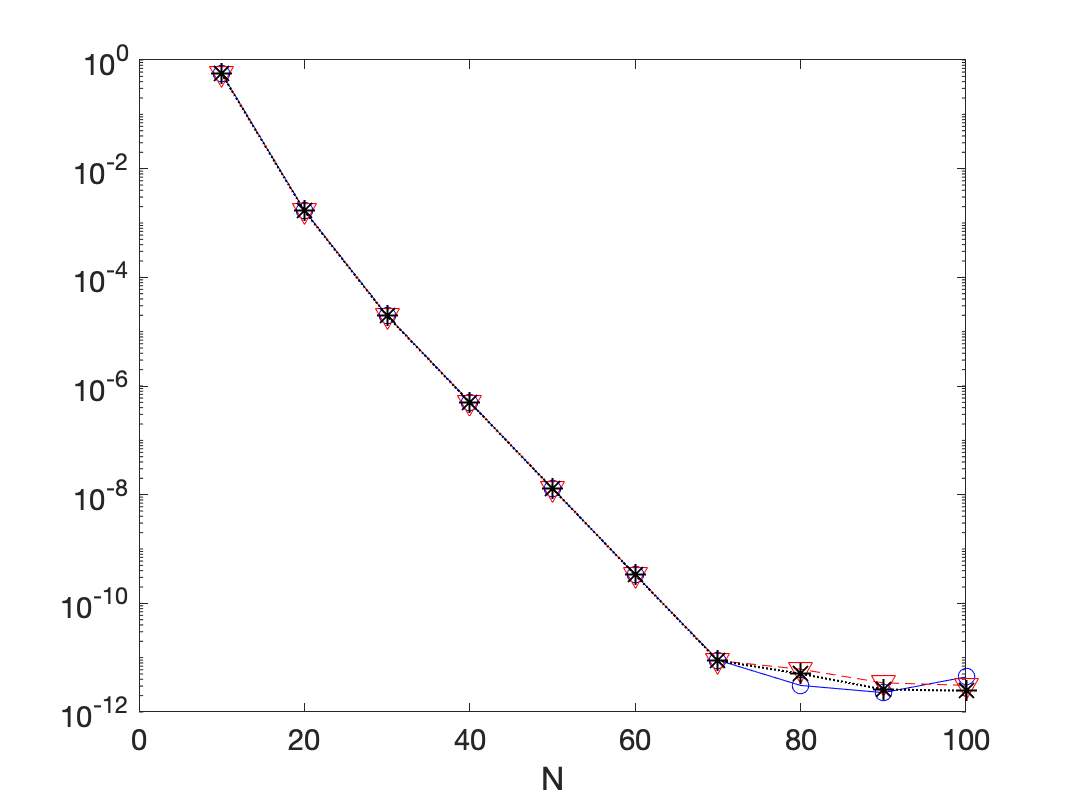}}
\caption{Maximum error (infinity norm) in the space semi-discretization w.r.t.~$N$, for the HBVM$(2s,s)$ methods used with time-step $h=0.1$ (left-plot) and $h=0.05$ (right-plot), for $s=1,2,3$. In both plots: the solid line with circles is used for HBVM(2,1); the dashed line with triangles is used for HBVM(4,2); the dotted line with stars is used for HBVM(6,3).}
\label{errN}
\end{figure}

\begin{figure}[t]
\centerline{\includegraphics[width=8cm,height=5.5cm]{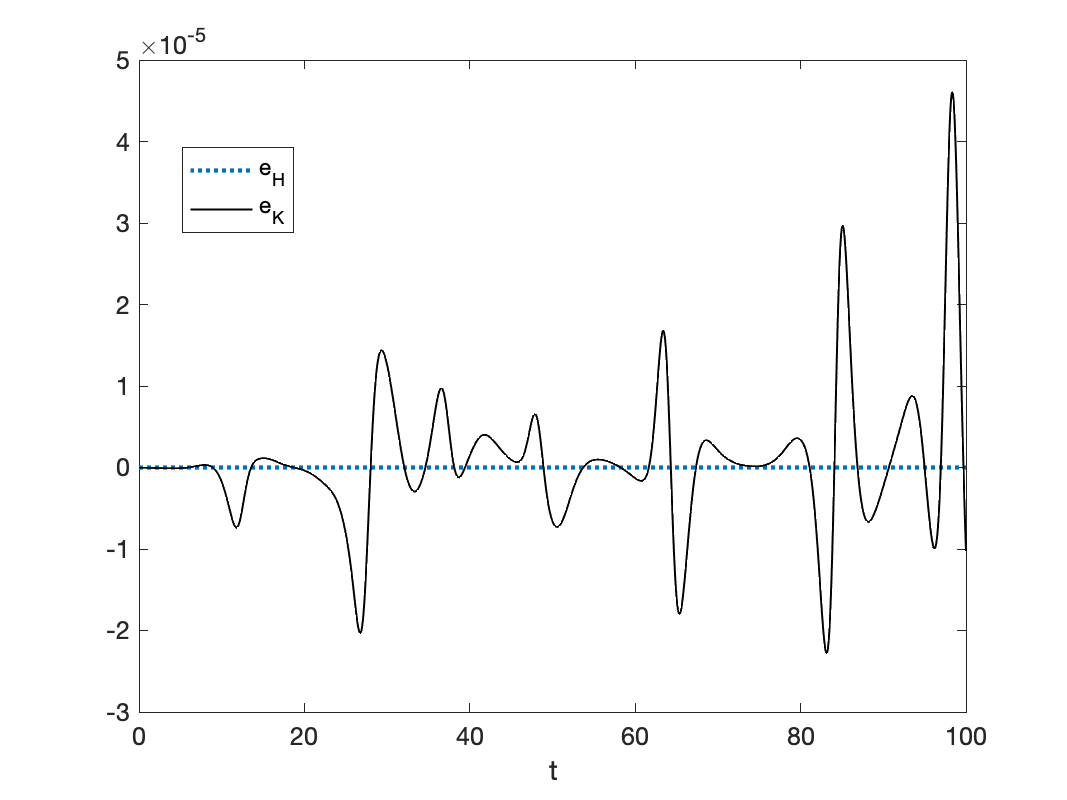}\includegraphics[width=8cm,height=5.5cm]{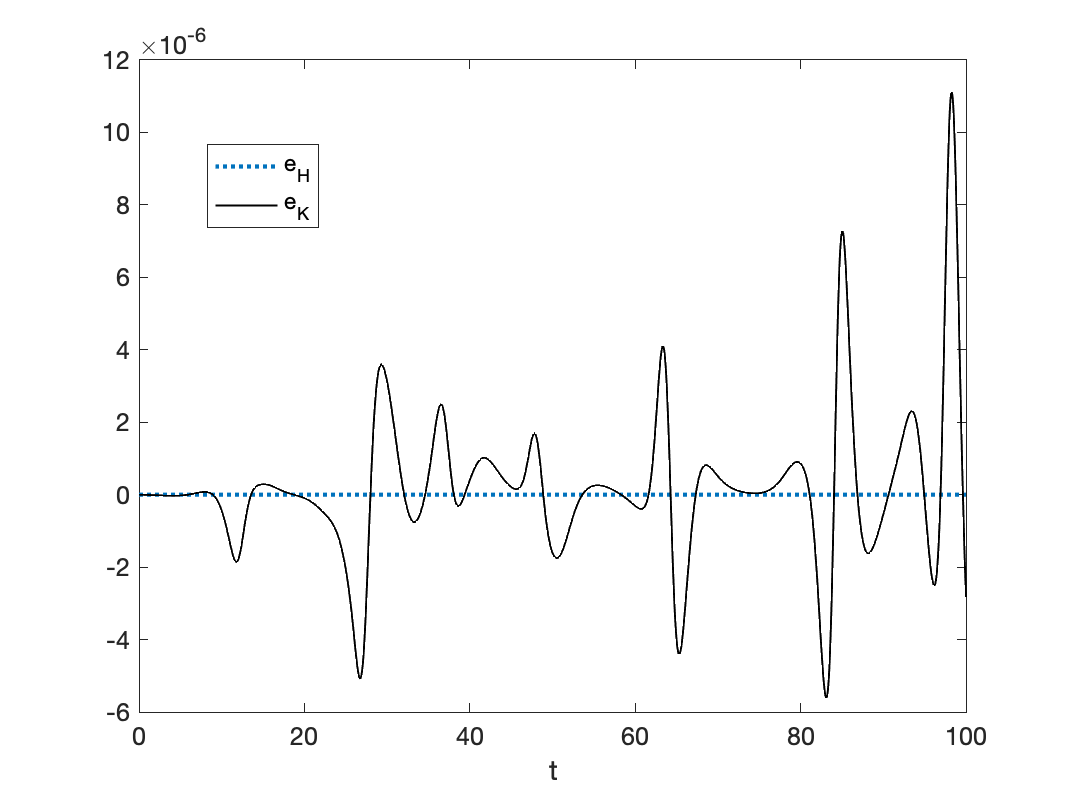}}
\centerline{\includegraphics[width=8cm,height=5.5cm]{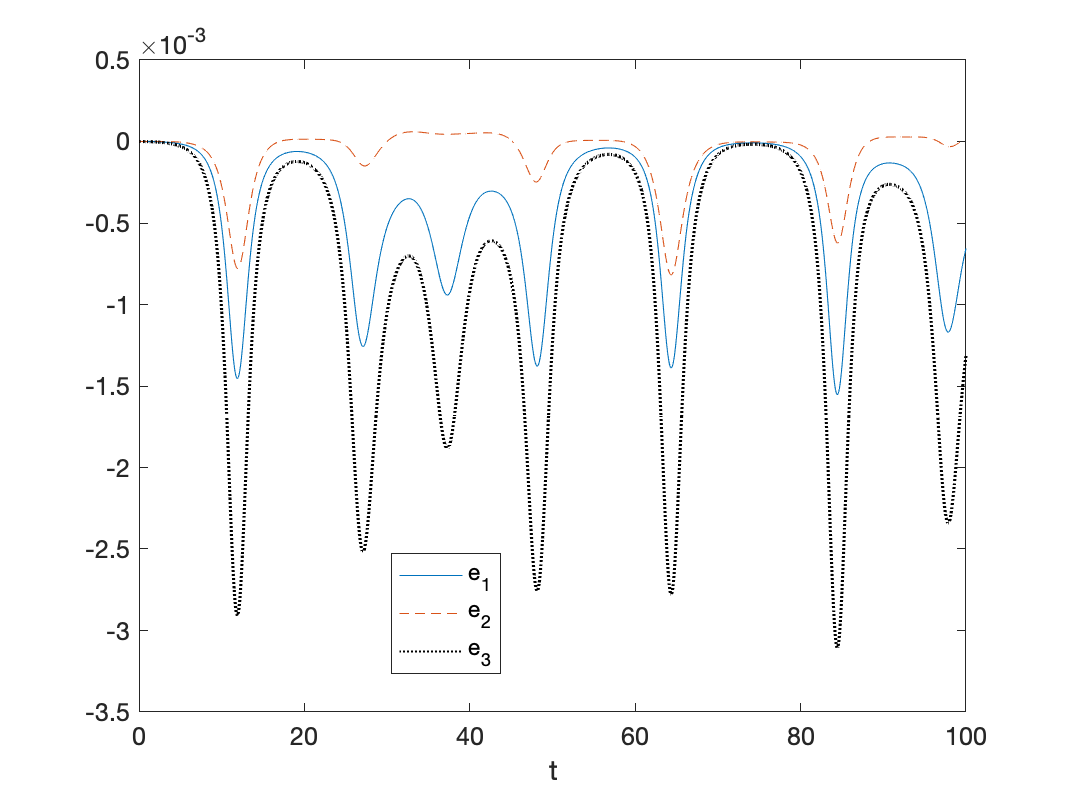}\includegraphics[width=8cm,height=5.5cm]{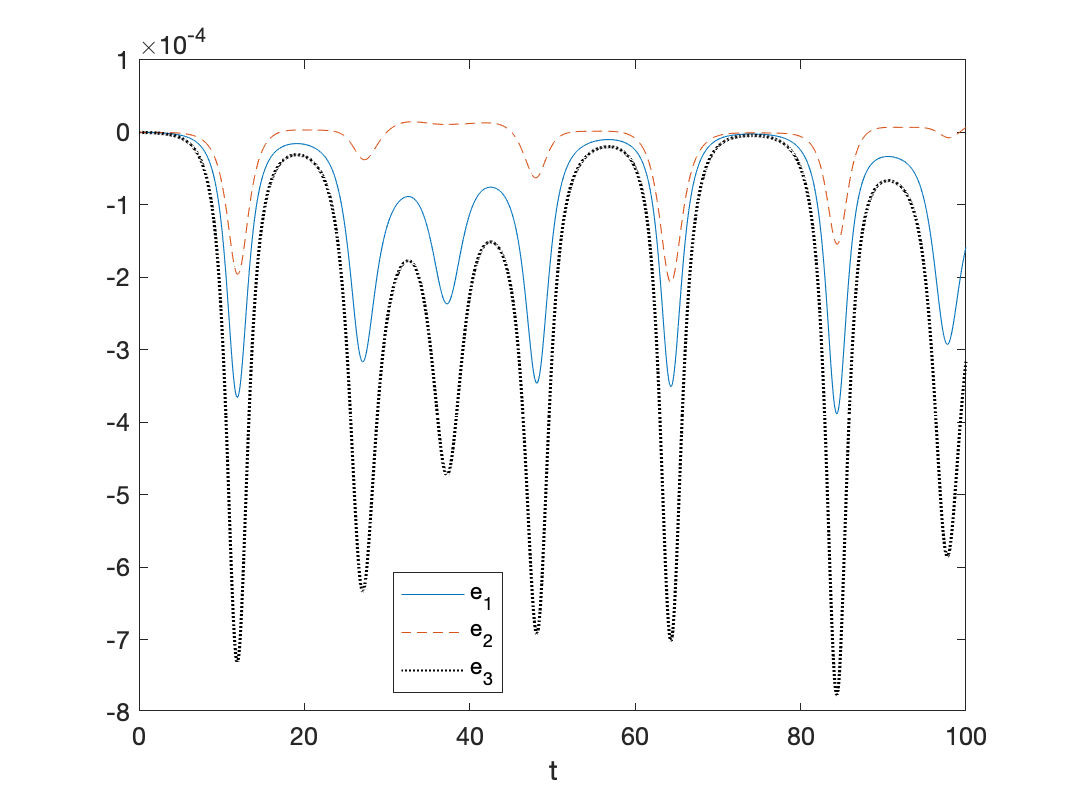}}
\centerline{\includegraphics[width=8cm,height=5.5cm]{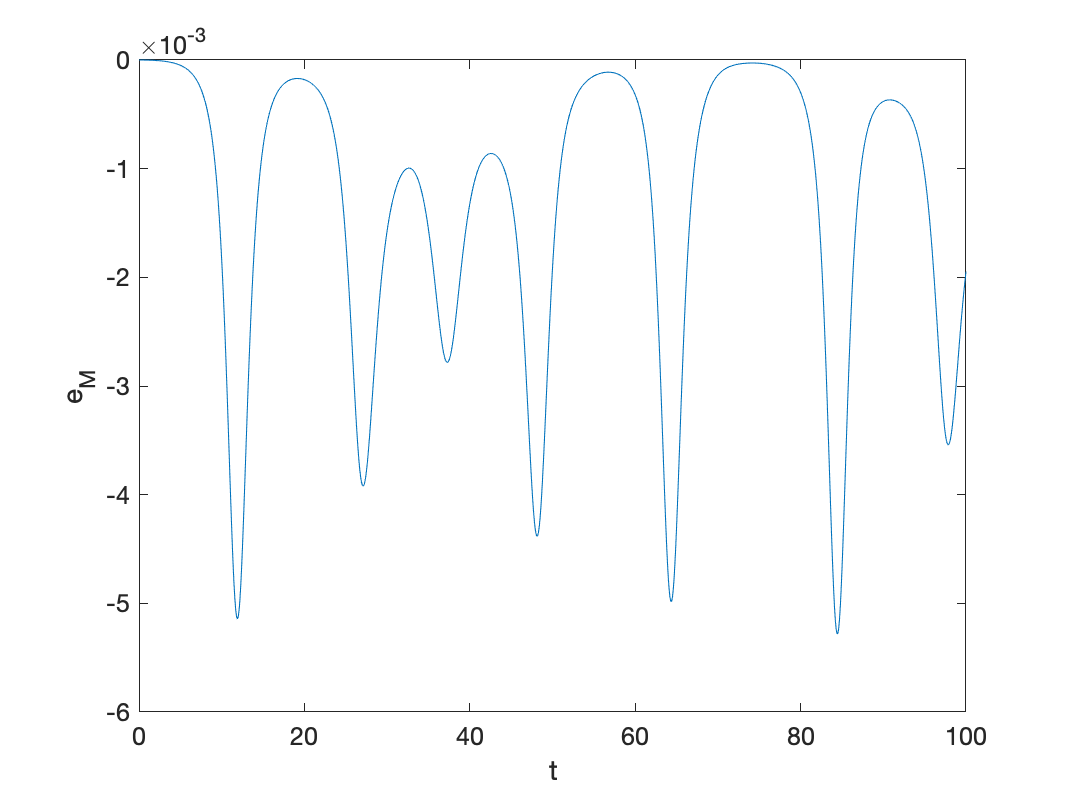}\includegraphics[width=8cm,height=5.5cm]{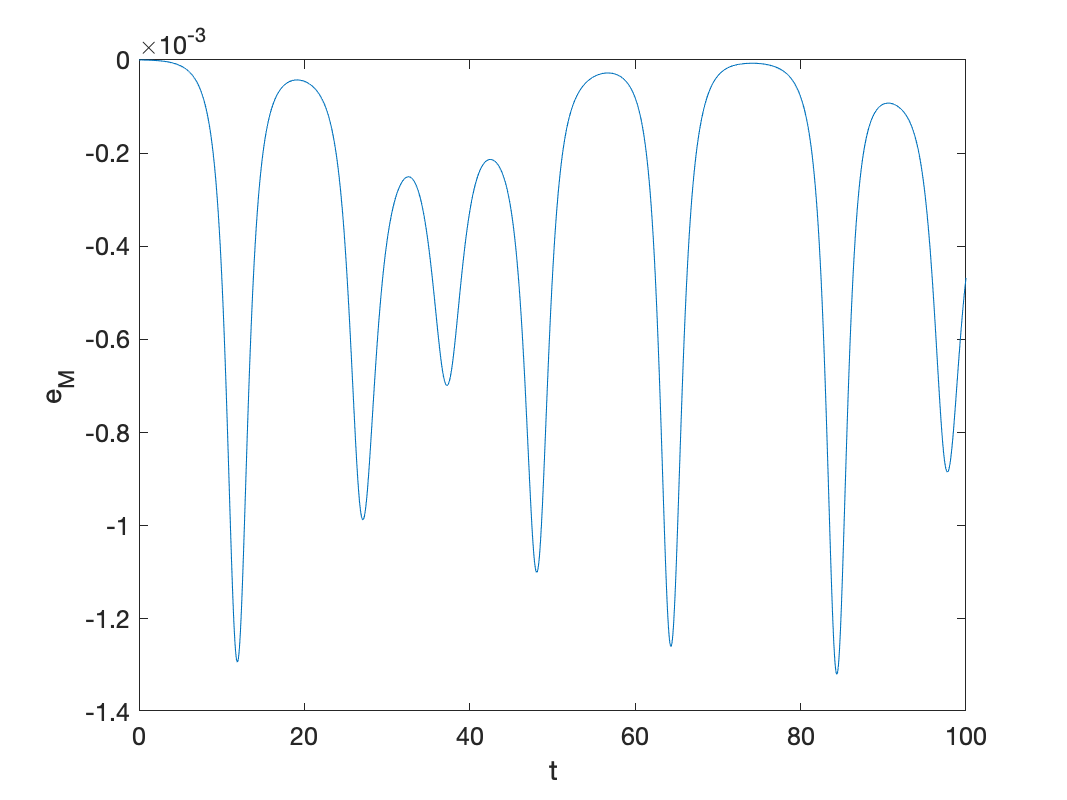}}
\caption{Invariant errors when solving problem (\ref{p1})-(\ref{p1par}) with the HBVM(2,1) method using time-steps $h=0.1$ (left-plots) and $h=0.05$ (right-plots). $e_H$ is the energy error; $e_K$ is the momentum error; $e_i$, $i=1,2,3$, is the $i$-th mass error; $e_M$ is the total mass error.}
\label{err21}
\end{figure}

\begin{figure}[t]
\centerline{\includegraphics[width=8cm,height=5.5cm]{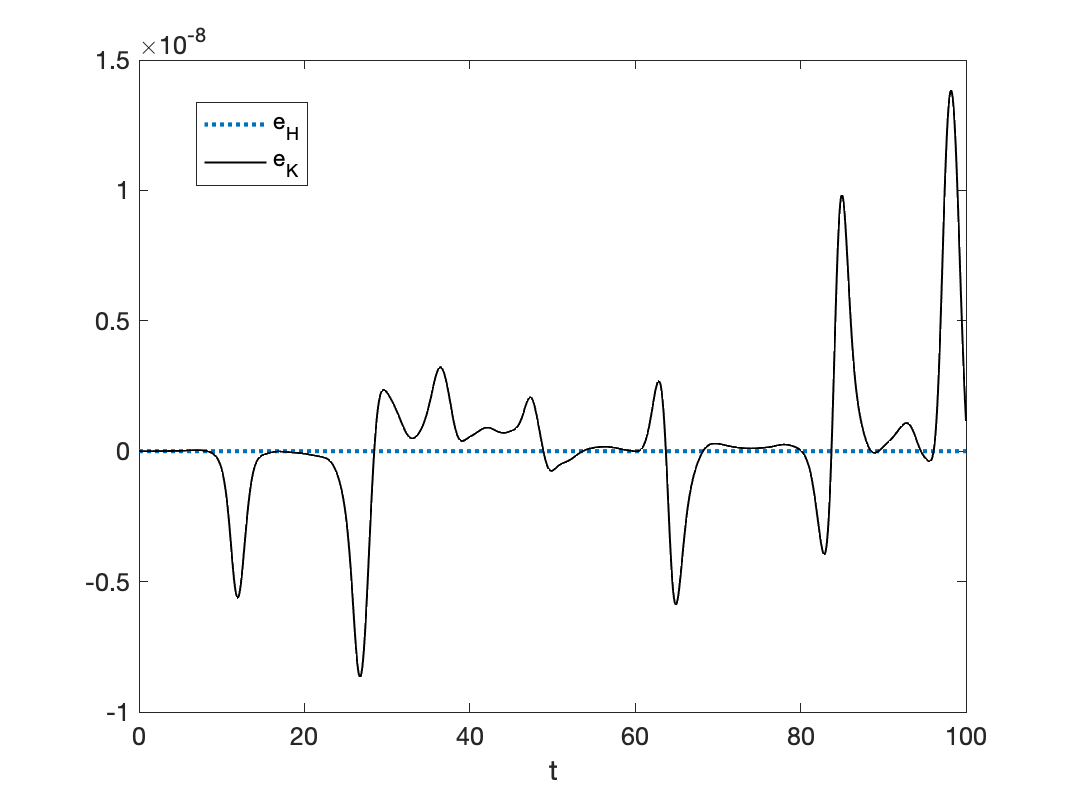}\includegraphics[width=8cm,height=5.5cm]{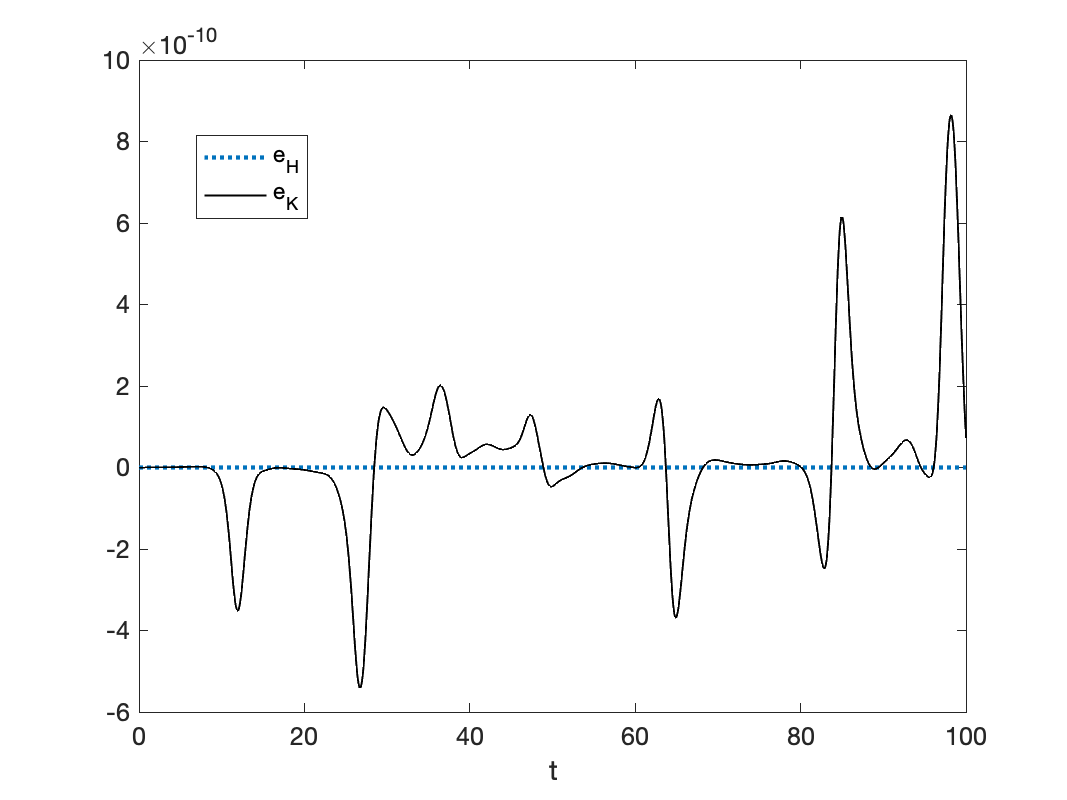}}
\centerline{\includegraphics[width=8cm,height=5.5cm]{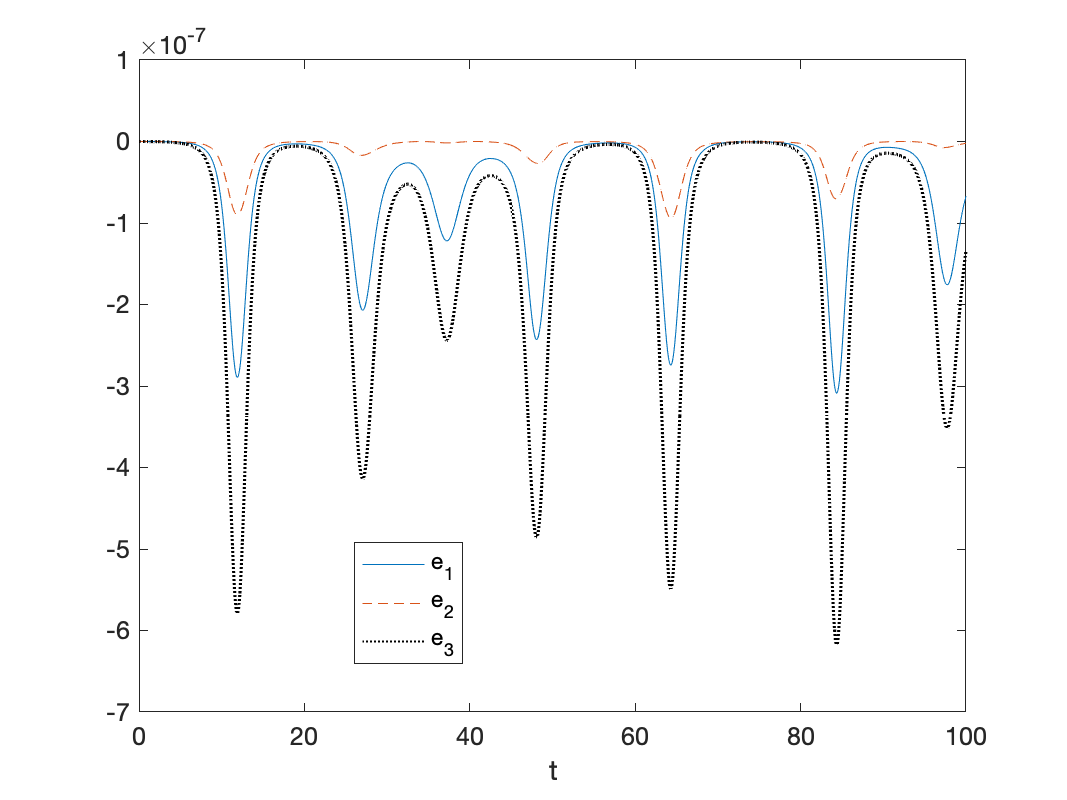}\includegraphics[width=8cm,height=5.5cm]{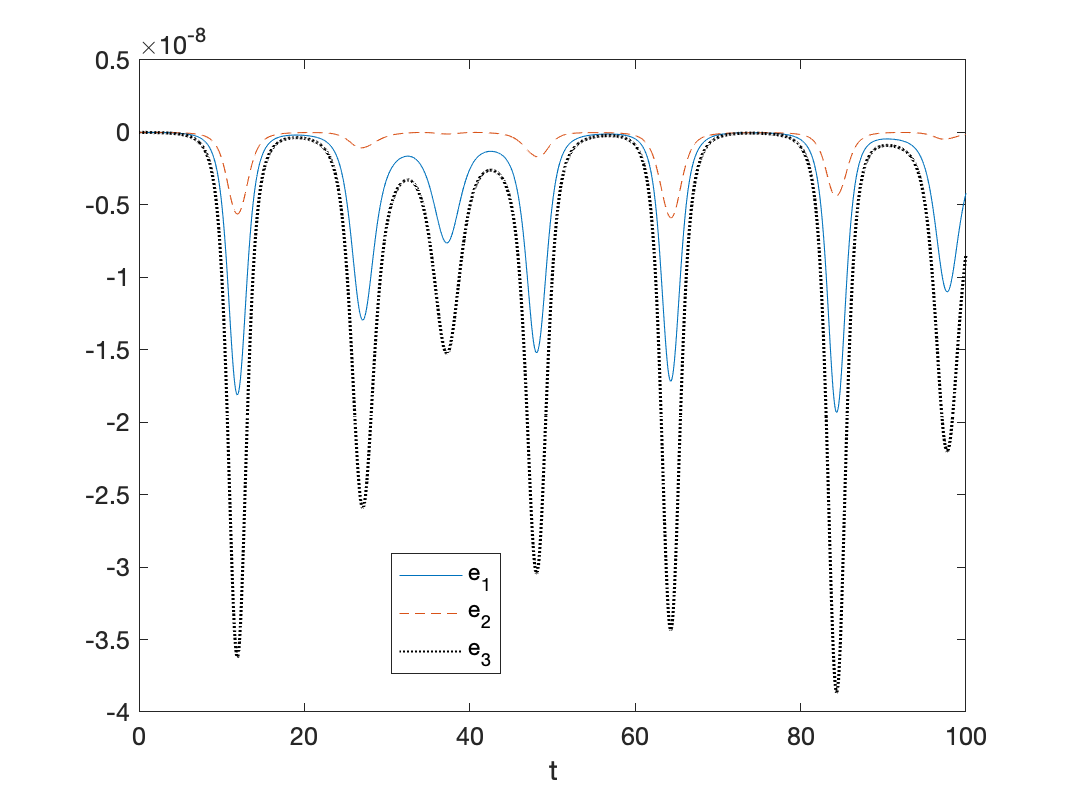}}
\centerline{\includegraphics[width=8cm,height=5.5cm]{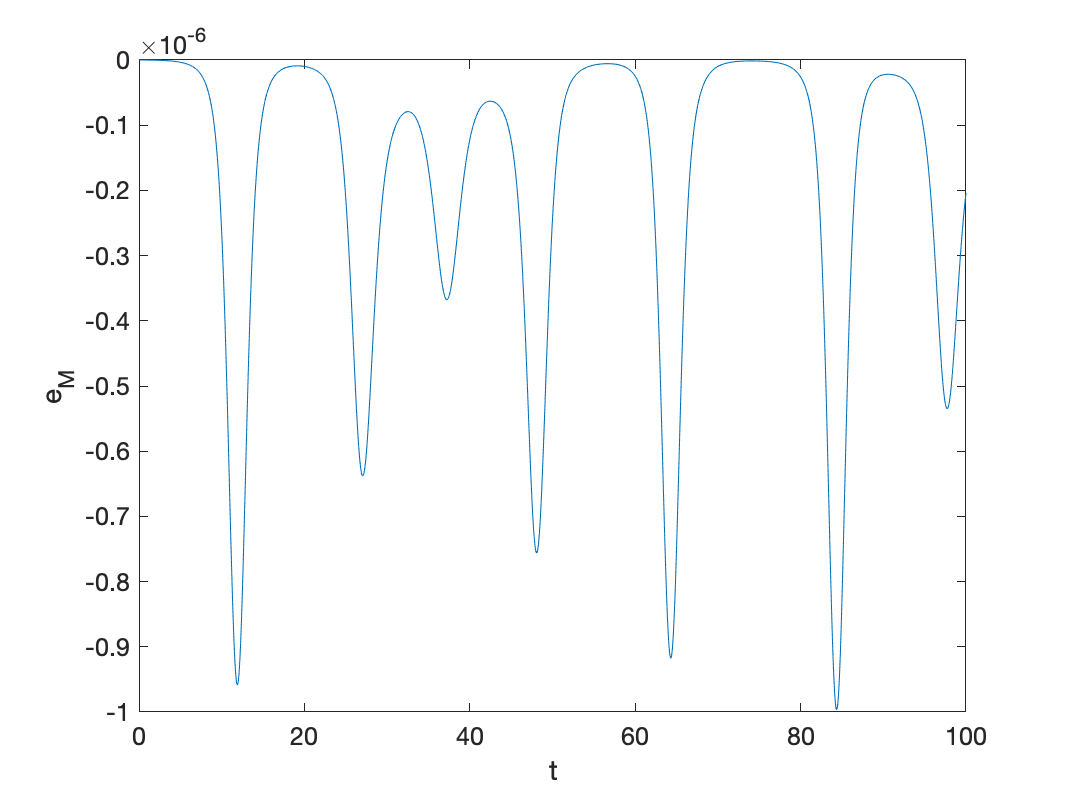}\includegraphics[width=8cm,height=5.5cm]{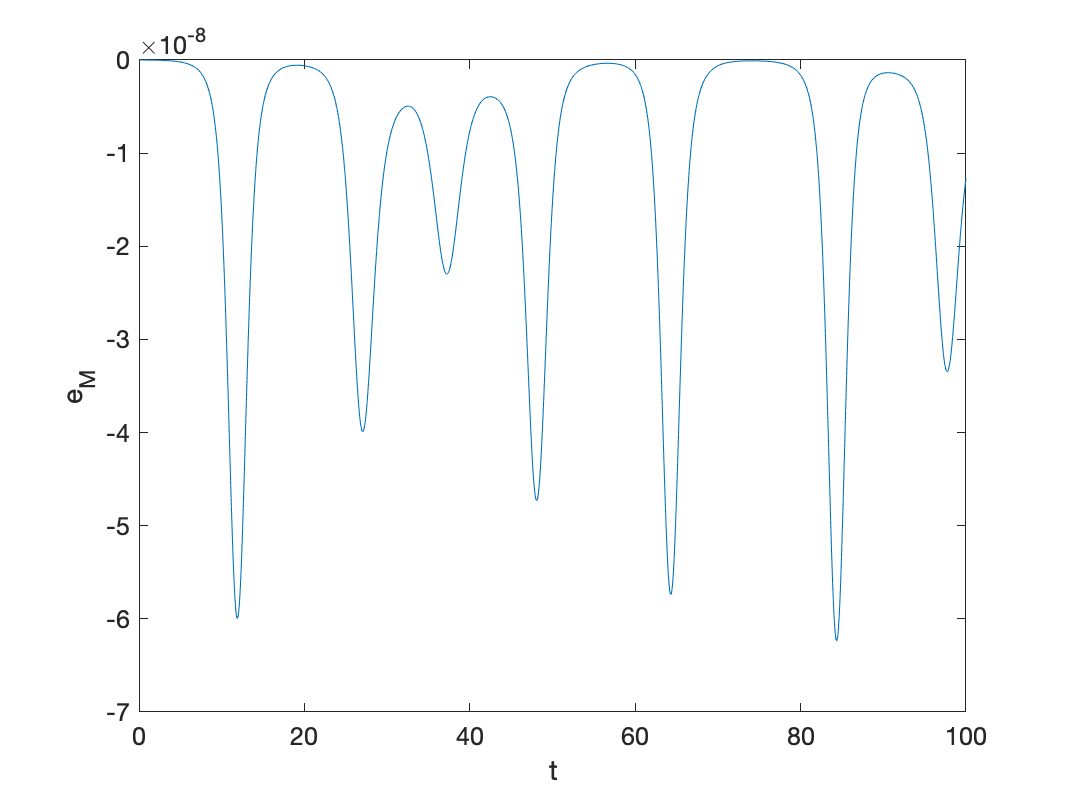}}
\caption{Invariant errors when solving problem (\ref{p1})-(\ref{p1par}) with the HBVM(4,2) method using time-steps $h=0.1$ (left-plots) and $h=0.05$ (right-plots). $e_H$ is the energy error; $e_K$ is the momentum error; $e_i$, $i=1,2,3$, is the $i$-th mass error; $e_M$ is the total mass error.}
\label{err42}
\end{figure}

\begin{figure}[t]
\centerline{\includegraphics[width=8cm,height=5.5cm]{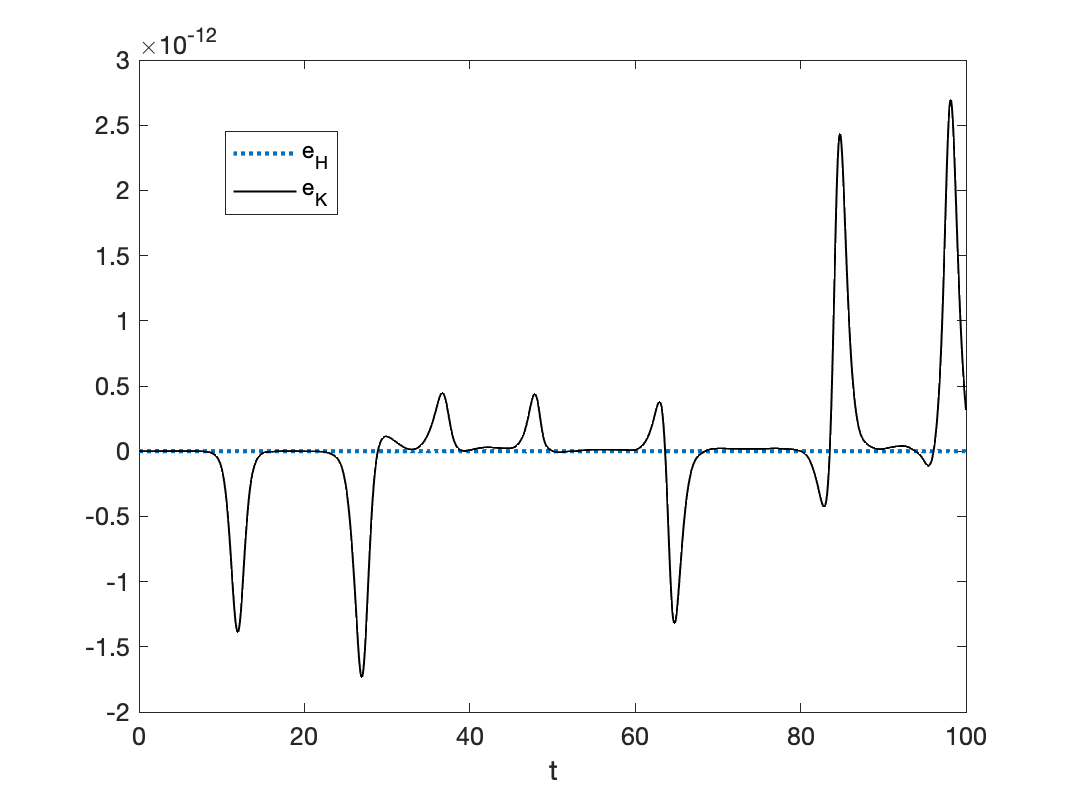}\includegraphics[width=8cm,height=5.5cm]{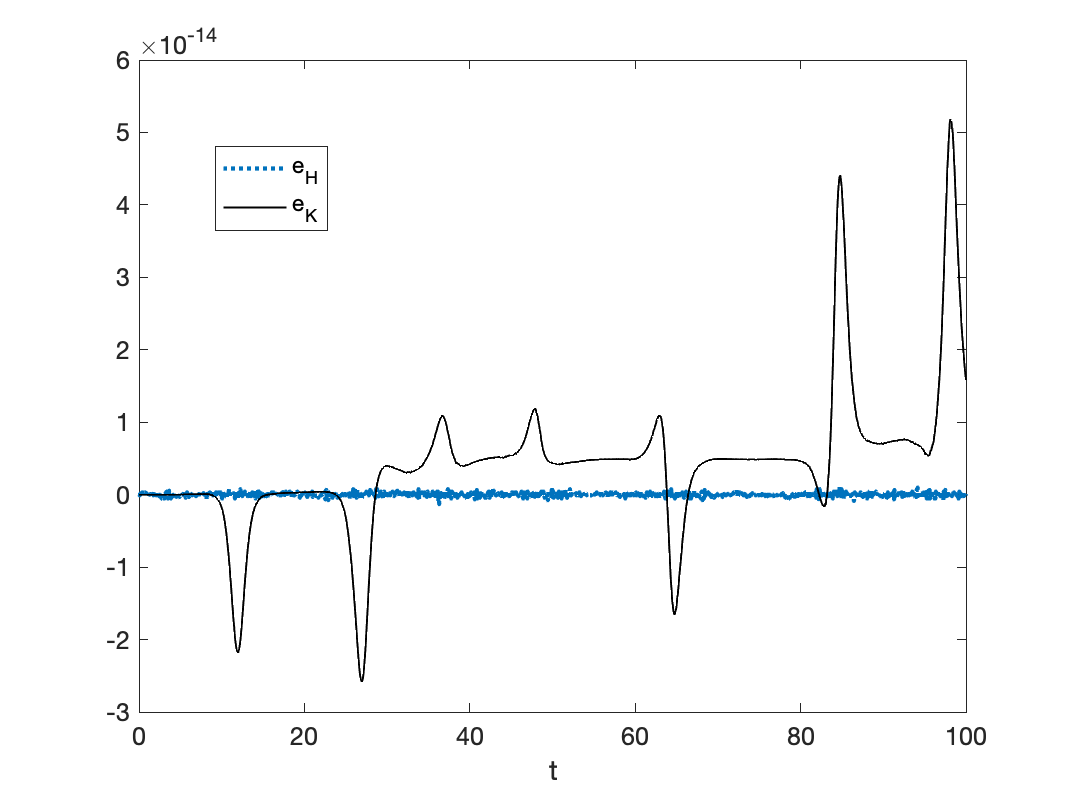}}
\centerline{\includegraphics[width=8cm,height=5.5cm]{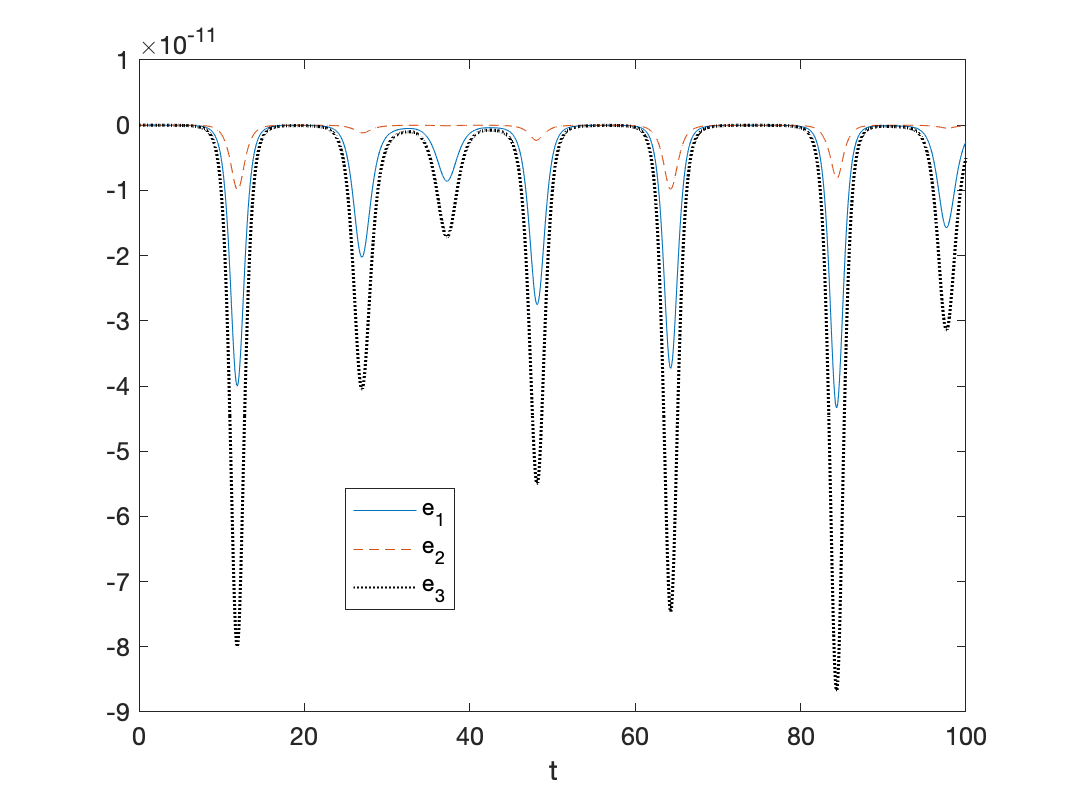}\includegraphics[width=8cm,height=5.5cm]{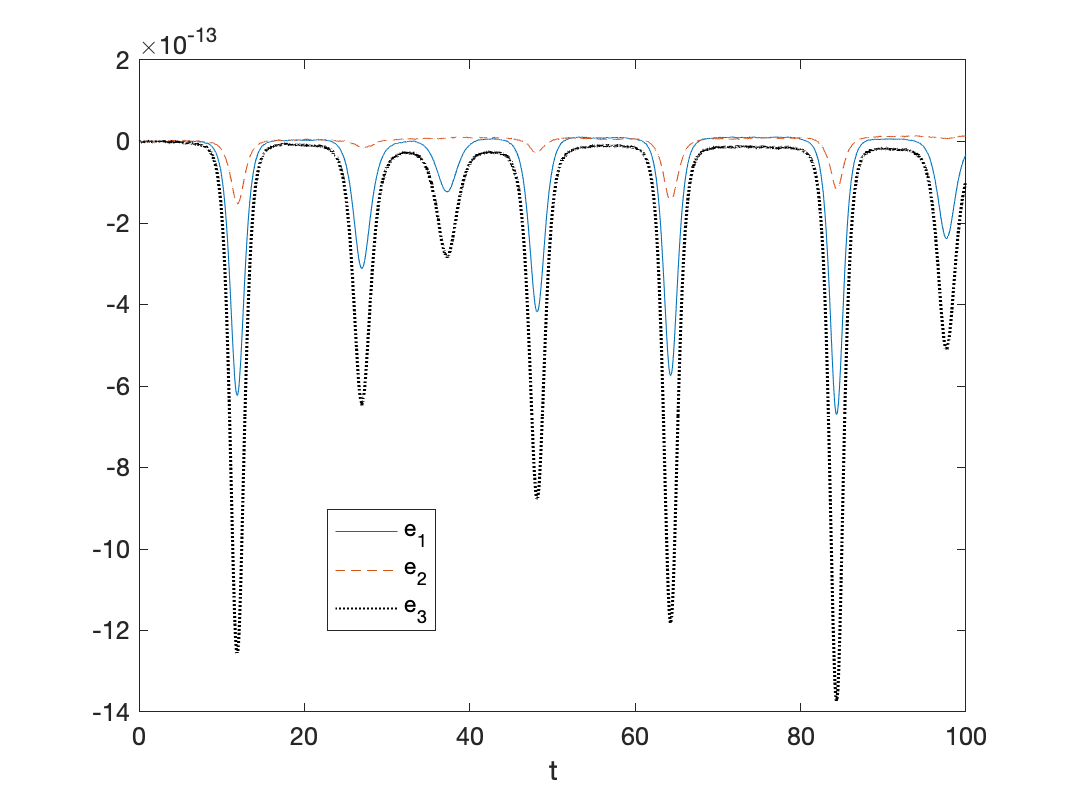}}
\centerline{\includegraphics[width=8cm,height=5.5cm]{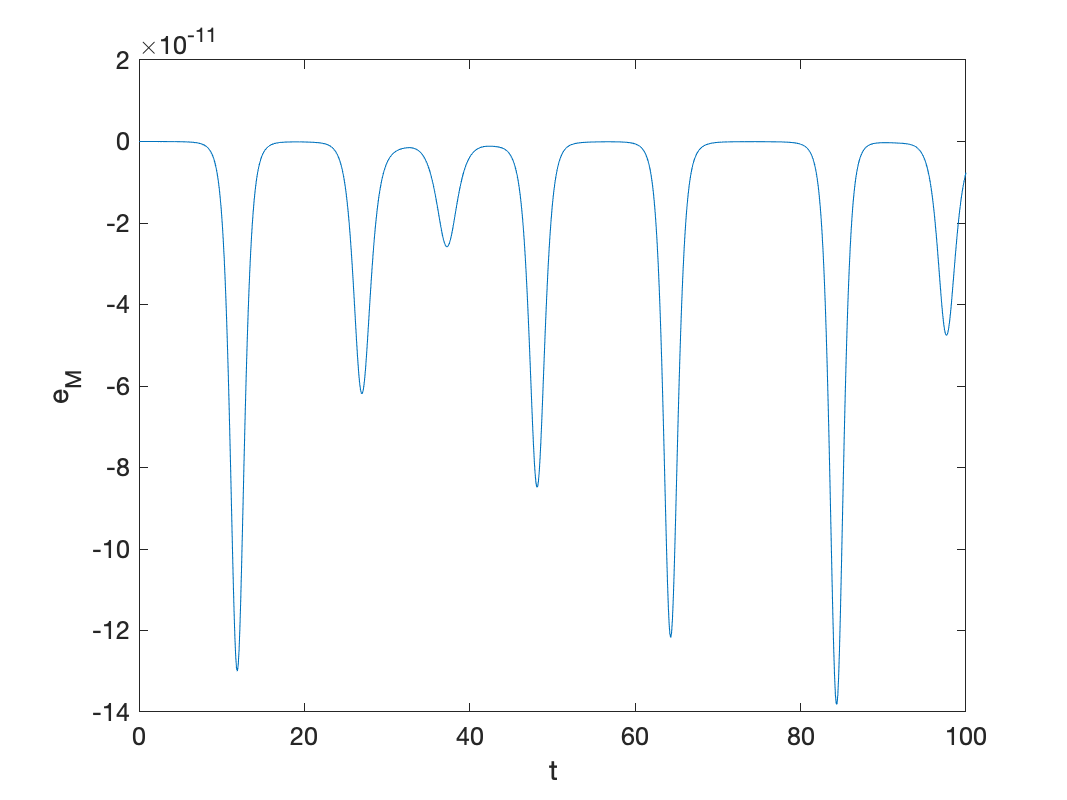}\includegraphics[width=8cm,height=5.5cm]{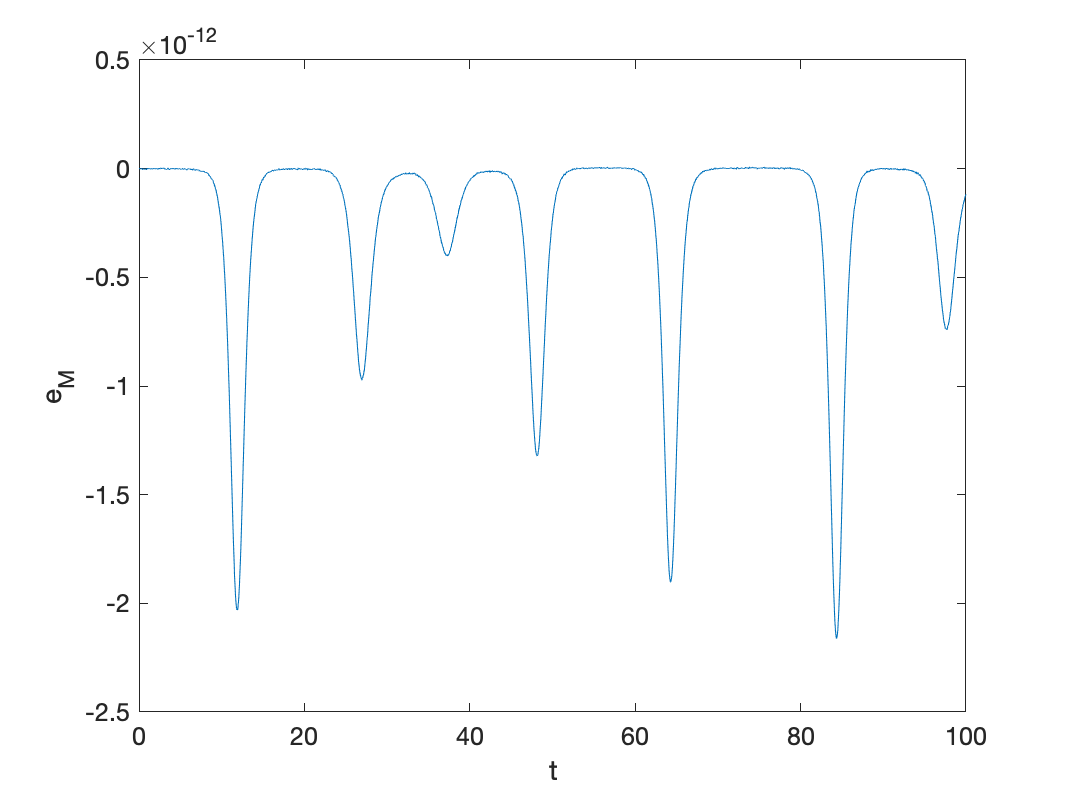}}
\caption{Invariant errors when solving problem (\ref{p1})-(\ref{p1par}) with the HBVM(6,3) method using time-steps $h=0.1$ (left-plots) and $h=0.05$ (right-plots). $e_H$ is the energy error; $e_K$ is the momentum error; $e_i$, $i=1,2,3$, is the $i$-th mass error; $e_M$ is the total mass error.}
\label{err63}
\end{figure}

The first test problem is taken from \cite{QSC2014}, and concern problems in the form (\ref{ms})-(\ref{ms1}), with
\begin{equation}\label{p1}
\beta = I_3, \qquad \gamma = \pmatrix{ccc} \sigma &e &\sigma\\ e&\sigma&e\\ \sigma & e &\sigma\endpmatrix,\qquad
\psi^0(x) = \pmatrix{c}
a_0(1-\eps\cos(\ell x))\\ 
b_0(1-\eps\cos(\ell( x+\theta)))\\
c_0(1-\eps\cos(\ell x))
\endpmatrix,\quad x\in[-4\pi,4\pi],
\end{equation}
and $T=100$. We here consider the following parameters, slightly different from those considered in the original reference, in order to emphasize some features:
\begin{equation}\label{p1par}
\sigma = 1, \quad e = 2/3\quad a_0 =  b_0 = 0.3, \quad c_0 = 0.3\sqrt{2}, \quad  \ell = 0.5,\quad   \eps = 0.1, \quad \theta = 9\pi/4.
\end{equation}
In fact, with such parameters the first two components have equal mass, which is one half of that of the third component. 
In Figure~\ref{solu} there is the plot of $|\psi_i(x,t)|^2$, $i=1,2,3$. The numerical tests for this problem are divided into four parts:
\begin{enumerate}
\item at first, for increasing values of $N$ used in the semi-discretization, we compute the error (in maximum norm) in the truncated solution over the whole integration interval by using the HBVM$(2s,s)$ methods with stepsizes $h=0.1$ and $0.05$, for $s=1,2,3$. The obtained results are reported in Figure~\ref{errN}, thus showing that a value $N=70$ is sufficient to obtain spectral accuracy in space.\footnote{The reference solution has been computed by using the given method with $N=150$.} This conclusion is further confirmed by the fact that the behaviour of the space semi-discretization error is the same for all methods and for all stepsizes used;

\item next, having fixed $N=70$ for the space semi-discretization, we use the HBVM$(2s,s)$ methods, $s=1,2,3$, with decreasing time-steps, and measure the maximum error (infinity norm), to assess the $2s$ order of accuracy (in time) of the methods. The obtained results are listed in Table~\ref{erroh}. As one may see, the expected order is confirmed;

\item next, in Figures~\ref{err21}--\ref{err63} we plot the errors in the invariants $H,K,M_1,M_2,M_3,M$ (see definitions (\ref{rMj1N})--(\ref{rH1N})), which we denote by $e_H,e_K,e_1,e_2,e_3,e_M$, respectively, for the HBVM(2,1), HBVM(4,2), HBVM(6,3) methods used with time-steps $h=0.1$ (left plots in each figure) and $h=0.05$ (right plots). As one may see, all methods are energy conserving.  Concerning the mass and momentum errors, one may see that no drift occurs and, moreover, they decrease with the expected order of the methods;

\item finally, according to the analysis in \cite{ABI2019}, in turn motivated by \cite{BIMR2018,BMR2018}, we show that the HBVM$(k,s)$ methods can be used as spectral methods in time, due to the availability of the very efficient blended iteration (\ref{blend}). For this purpose, in Table~\ref{spectral} we list the maximum solution error, $e_y$, the Hamiltonian error, $e_H$, the momentum error, $e_K$, and the total mass error, $e_M$, for various combinations of the order and of the time-steps used to cover the integration interval $[0,100]$ (as before, the value $N=70$ has been fixed for the space semi-discretization). HBVM(20,10) amounts to the use of the method (for the given time-step $h=1$) as a spectral method in time, according to what sketched in Section~\ref{spectrem}.   As one may see, this is the most efficient one, allowing the use of very large time-step ($h=1$, in the present case), and a very small execution time.  This further confirms what observed in \cite{BIMR2018}. 

\end{enumerate}

\begin{table}[t]
\caption{Order of convergence in time for the HBVM$(2s,s)$ method solving problem (\ref{p1})-(\ref{p1par}) with $N=70$ and time-steps $h_i=2^{1-i}/10$; $e_y$ denotes the maximum solution error.}\label{erroh} 

\smallskip
\centerline{\begin{tabular}{|c|rr|rr|rr|}
\hline
\hline
      &\multicolumn{2}{|c|}{HBVM(2,1)} &\multicolumn{2}{|c|}{HBVM(4,2)} &\multicolumn{2}{|c|}{HBVM(6,3)}\\
\hline      
$i$ & $e_y$~~ & rate & $e_y$~~ & rate & $e_y$~~ & rate \\
\hline
   0&   3.712e-01&    ---&    2.877e-04&   ---&    2.646e-07&   ---\\
   1&   1.055e-01&   1.8&   1.814e-05&   4.0&   4.108e-09&   6.0\\
   2&   2.715e-02&   2.0&   1.135e-06&   4.0&   6.399e-11&   6.0\\
   3&   6.833e-03&   2.0&   7.099e-08&   4.0&   1.023e-12&   6.0\\
\hline
\hline
\end{tabular}
}
\end{table}


\begin{table}[t]
\caption{Numerical solution of problem (\ref{p1})-(\ref{p1par}) by using HBVMs with time-step $h$: $e_y$ is the solution error (maximum norm); $e_H$, $e_K$, and $e_M$ are the errors in the invariants. The execution times are in {\em sec}.}
\label{spectral} 
\smallskip
\centerline{\begin{tabular}{|rrrrrrrrr|}
\hline
\hline
\multicolumn{9}{|c|}{HBVM(2,1)} \\
\hline
$h$          & $e_y$ & rate & $e_H$ & $e_K$ & rate & $e_M$ & rate & time \\
\hline
1.000e-01 & 1.055e-01 & --- & 1.332e-15 & 4.604e-05 & --- & 5.280e-03 & --- & 11.7 \\ 
5.000e-02 & 2.715e-02 & 2.0 & 1.332e-15 & 1.111e-05 & 2.0 & 1.319e-03 & 2.0 & 17.4 \\ 
2.500e-02 & 6.833e-03 & 2.0 & 1.665e-15 & 2.753e-06 & 2.0 & 3.296e-04 & 2.0 & 31.3 \\ 
1.250e-02 & 1.711e-03 & 2.0 & 1.554e-15 & 6.866e-07 & 2.0 & 8.239e-05 & 2.0 & 57.9 \\ 
6.250e-03 & 4.279e-04 & 2.0 & 1.554e-15 & 1.716e-07 & 2.0 & 2.060e-05 & 2.0 & 105.9 \\ 
3.125e-03 & 1.070e-04 & 2.0 & 1.665e-15 & 4.288e-08 & 2.0 & 5.149e-06 & 2.0 & 181.1 \\ 
\hline
\multicolumn{9}{|c|}{HBVM(4,2)} \\
\hline
$h$          & $e_y$ & rate & $e_H$ & $e_K$ & rate & $e_M$ & rate & time \\
\hline
1.000e-01 & 1.814e-05 & --- & 1.554e-15 & 1.383e-08 & --- & 9.962e-07 & --- & 11.9 \\ 
5.000e-02 & 1.135e-06 & 4.0 & 1.554e-15 & 8.647e-10 & 4.0 & 6.236e-08 & 4.0 & 20.3 \\ 
2.500e-02 & 7.099e-08 & 4.0 & 1.332e-15 & 5.401e-11& 4.0 & 3.898e-09 & 4.0 & 36.7 \\ 
1.250e-02 & 4.437e-09 & 4.0 & 1.665e-15 & 3.336e-12 & 4.0 & 2.436e-10 & 4.0 & 64.4 \\ 
6.250e-03 & 2.774e-10 & 4.0 & 1.887e-15 & 1.656e-13 & 4.3 & 1.521e-11 & 4.0 & 201.0 \\ 
\hline
\multicolumn{9}{|c|}{HBVM(6,3)} \\
\hline
$h$          & $e_y$ & rate & $e_H$ & $e_K$ & rate & $e_M$ & rate & time \\
\hline
1.000e-01 & 4.108e-09 & --- & 1.110e-15 & 2.640e-12 & --- & 1.381e-10 & --- & 18.3 \\ 
5.000e-02 & 6.399e-11 & 6.0 & 1.332e-15 & 5.279e-14 & 5.6  & 2.162e-12 & 6.0 & 29.0 \\ 
2.500e-02 & 1.023e-12 & 6.0 & 1.554e-15 & 4.371e-14 & *** & 3.730e-14 & 5.9 & 82.4 \\ 
\hline
\multicolumn{8}{|c|}{HBVM(20,10)} \\
\hline
$h$          & $e_y$ &  --- & $e_H$ & $e_K$ & --- & $e_M$ & --- & time \\
\hline
1.0           & 6.365e-11 &  & 1.332e-15 & 1.127e-14 & & 1.066e-14 & & 13.4\\
\hline
\hline
\end{tabular}
}
\end{table}

\subsection*{Second test problem}  
The second set of test problems, consisting in solitary waves (or vector solitons), is adapted from \cite{Ismail2008}, and is defined by:
\begin{equation}\label{p2}
\beta = \frac{1}2I_3, \qquad \gamma = (1+e)\pmatrix{ccc} 1& 1& 1\\ 1 &1 &1\\ 1& 1&1\endpmatrix,\qquad T=40,
\end{equation}
with the initial condition
\begin{equation}\label{p2ini}
\psi^0(x) = \pmatrix{c}
\sqrt{\frac{2\alpha_1}{1+e}}\,\sech\left(\sqrt{2\alpha_1}(x-x_1)\right)\exp\left(\ii v_1(x-x_1)\right)\\
\sqrt{\frac{2\alpha_2}{1+e}}\,\sech\left(\sqrt{2\alpha_2}(x-x_2)\right)\exp\left(\ii v_2(x-x_2)\right)\\
\sqrt{\frac{2\alpha_3}{1+e}}\,\sech\left(\sqrt{2\alpha_3}(x-x_3)\right)\exp\left(\ii v_3(x-x_3)\right)
\endpmatrix,\qquad x\in[-20,85].
\end{equation}
We use the parameters
\begin{eqnarray}\label{p2par}
&&e = \frac{2}3, \qquad \alpha_1 =  1,\qquad \alpha_2 = 0.6, \qquad \alpha_3 = 0.3, \\[2mm] \nonumber
&&v_1 = 1,\quad v_2 = 0.1,\quad v_3 = -1,\quad x_1 = 0,\quad x_2 = 22,\quad x_3 = 50.
\end{eqnarray}

\begin{figure}[t]
\centerline{\includegraphics[width=8cm,height=6.5cm]{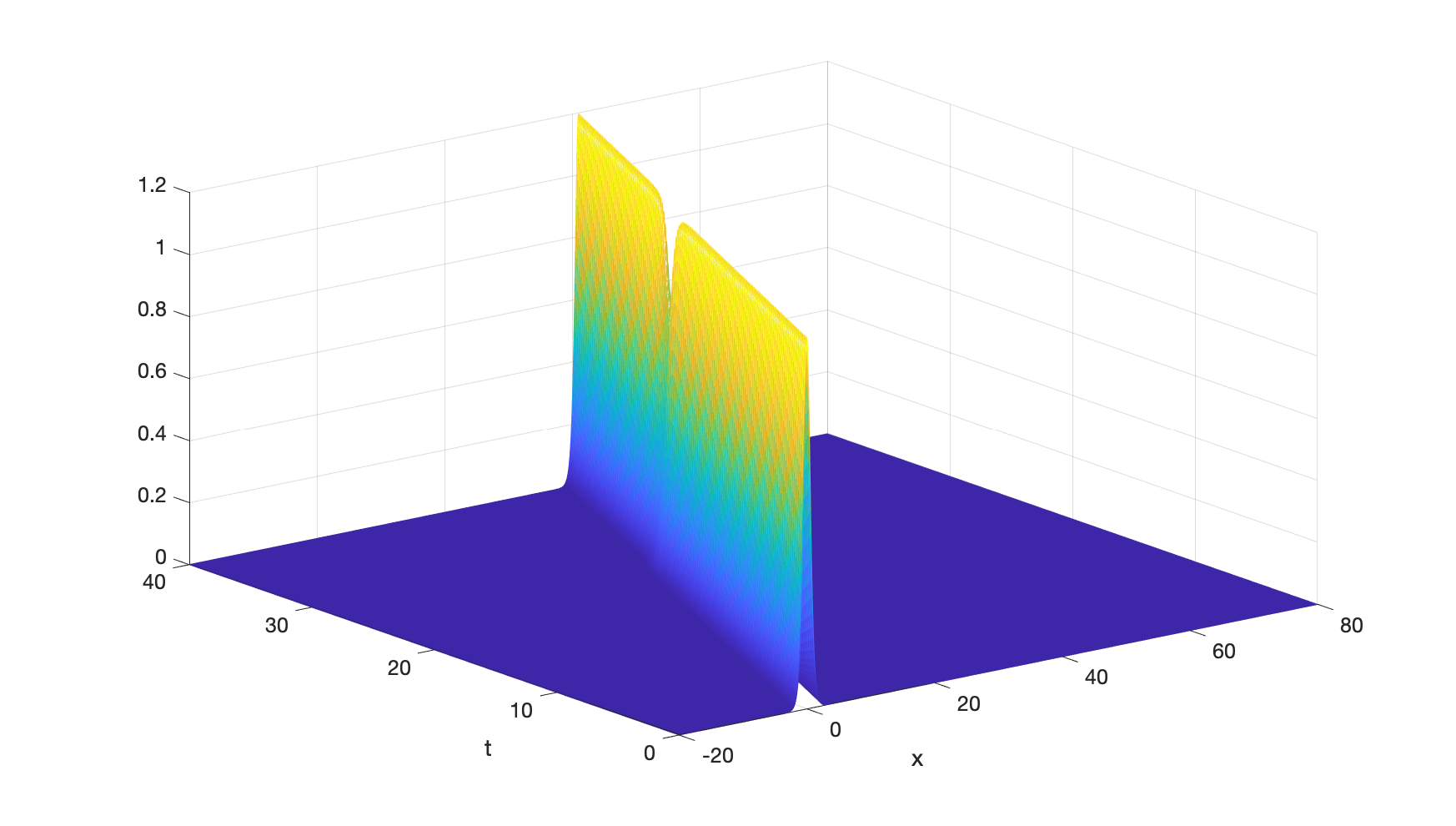}\includegraphics[width=8cm,height=6.5cm]{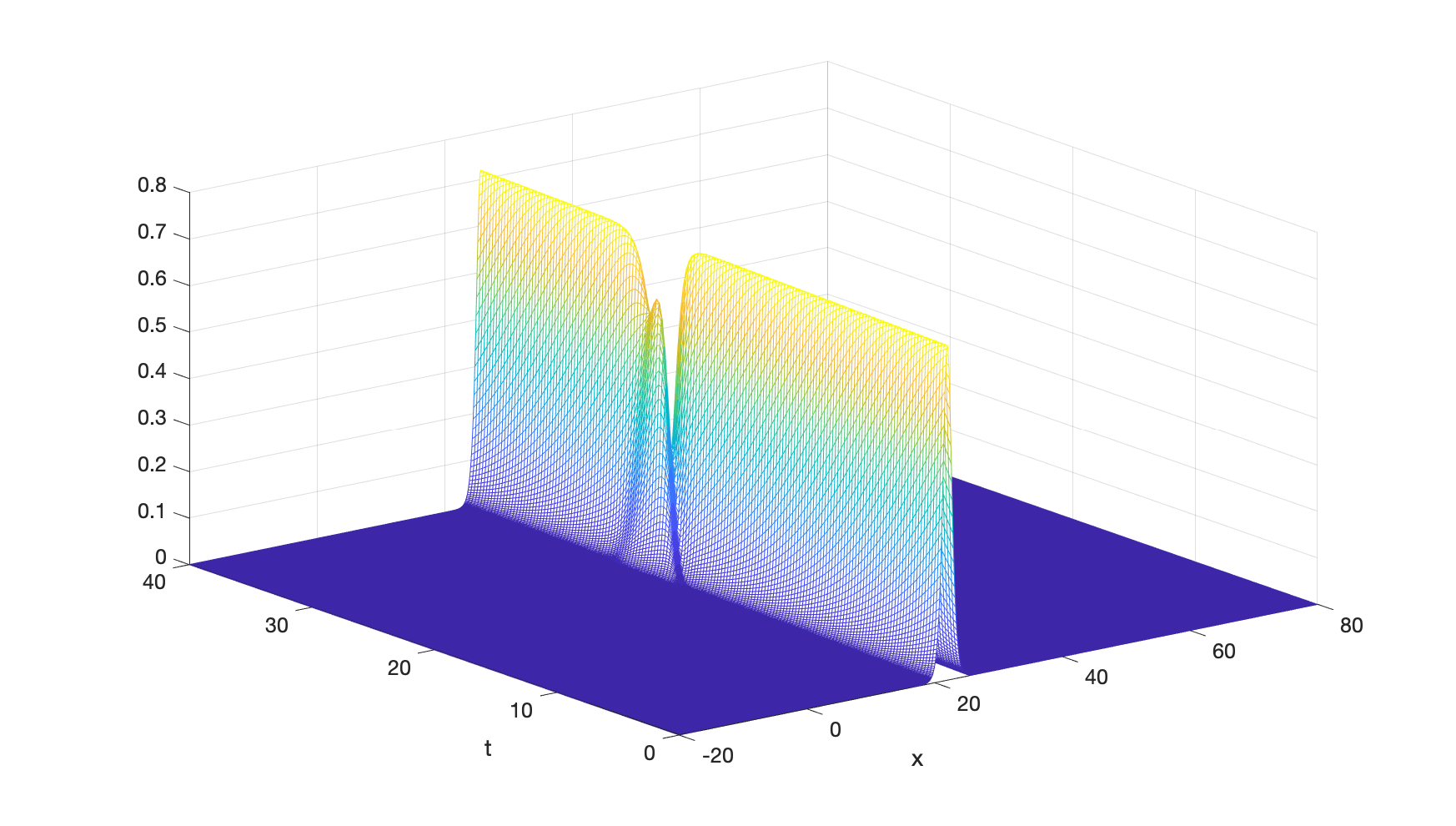}}
\centerline{\includegraphics[width=8cm,height=6.5cm]{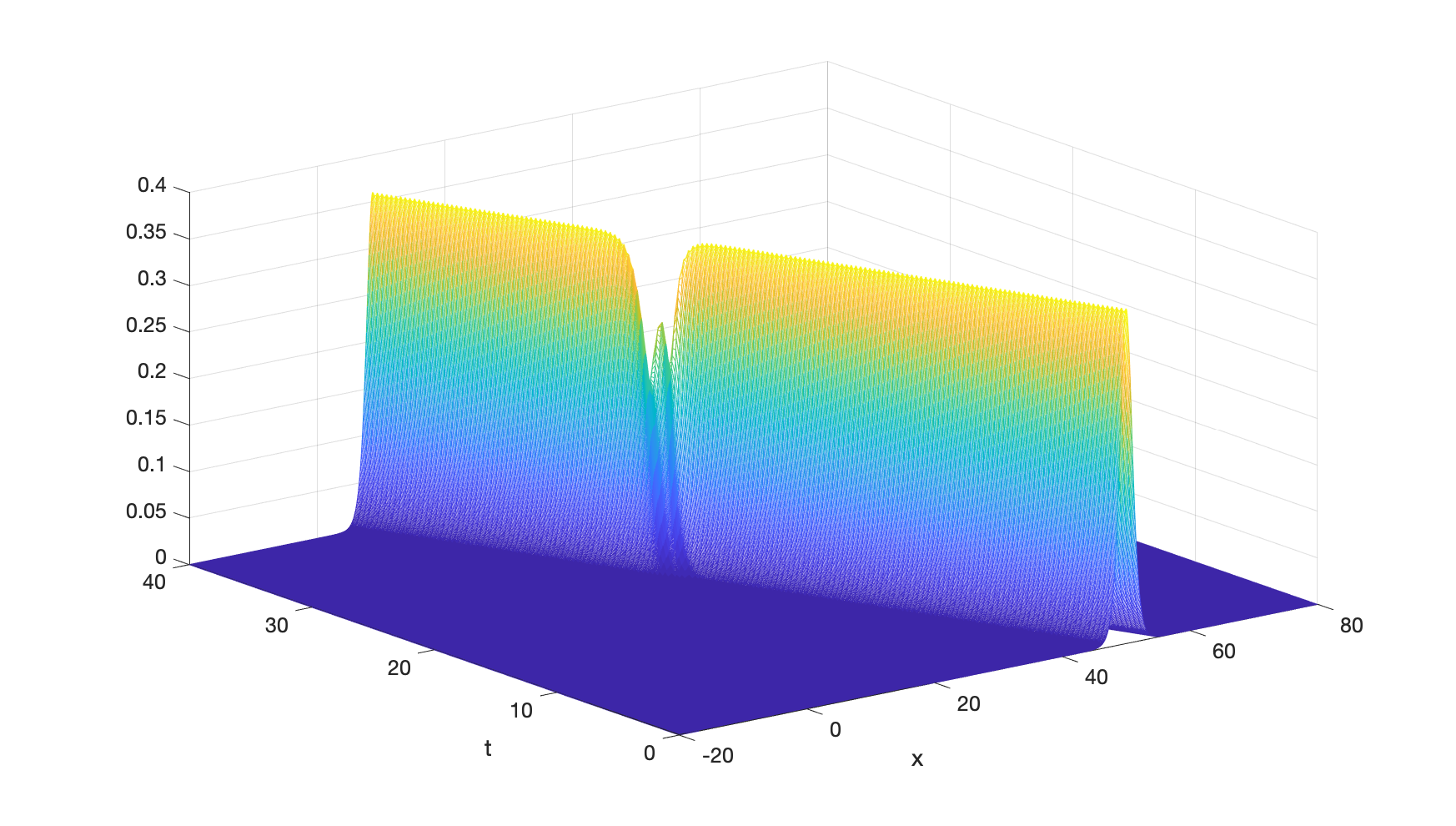}\includegraphics[width=8cm,height=6.5cm]{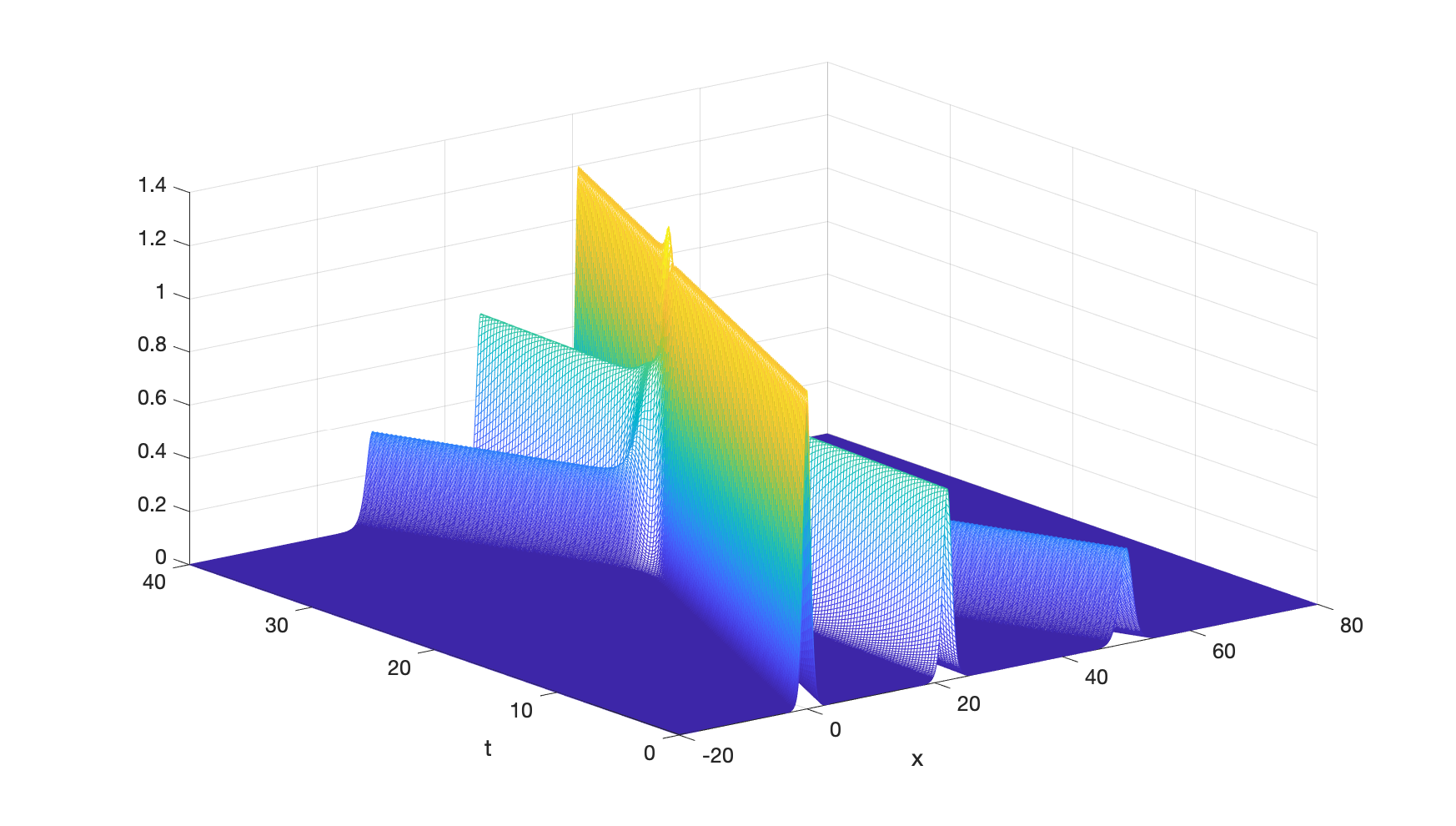}}
\caption{Left to right and up to down, plots of $|\psi_i(x,t)|^2$, $i=1,2,3$, and of $\sum_i|\psi_i(x,t)|^2$, for problem (\ref{p2})--(\ref{p2par}).}
\label{solup2}
\end{figure}

In Figure~\ref{solup2} is the plot of $|\psi_i(x,t)|^2$, $i=1,2,3$, and of their sum, for $(x,t)\in[-20,80]\times[0,40]$: as one may see, the components of the vector soliton are almost independent, with the exception of a small region near $x\approx t\approx 24$, where they have a significant interaction, due to the coupling nonlinear term, after which the solitons re-emerge unchanged. In this case, a value $N=400$ is required to obtain a spectrally accurate space semi-discretization (with an error of the order of  $6\cdot 10^{-12}$). In Table~\ref{spectral1} we list the obtained results by using the following methods:
\begin{itemize}
\item the symplectic $s$-stage Gauss method (i.e., HBVM$(s,s)$), $s=1,2$, which exactly conserves all quadratic invariants \cite{SS1988};
\item the energy-conserving HBVM$(2s,s)$ methods, $s=1,2$;
\item the HBVM(20,16) method, which provides a spectrally accurate time solution, for the given stepsize ($h=1$, in the present case).
\end{itemize}
\begin{table}[t]
\caption{Numerical solution of problem (\ref{p2})--(\ref{p2par}) by using HBVM$(k,s)$ with time-step $h$: $e_y$ is the solution error (maximum norm); $e_H$, $e_K$, and $e_M$ are the errors in the invariants, blend is the total number of {\em blended iterations} (\ref{blend}). The execution times are in {\em sec}.}
\label{spectral1} 
\smallskip
\centerline{\begin{tabular}{|rrrrrrrrrrr|}
\hline
\hline
\multicolumn{11}{|c|}{HBVM(1,1)} \\
\hline
$h$          & $e_y$ & rate & $e_H$ & rate & $e_K$ & --- & $e_M$ & --- & blend & time \\
\hline
1.00e-01 & 2.755e-01 & ---  & 1.143e-03 & ---   &7.234e-13 &&  5.249e-13 & & 6014 & 22.3 \\ 
5.00e-02 & 7.277e-02 & 1.9 & 2.903e-04 & 2.0 &7.105e-13 &&  3.002e-13   & &8846 & 48.3 \\ 
2.50e-02 & 1.845e-02 & 2.0 & 7.291e-05 & 2.0 &7.272e-14 &&  4.086e-14   & &14400 & 81.5 \\ 
1.25e-02 & 4.628e-03 & 2.0 & 1.825e-05 & 2.0 &1.477e-14 &&  1.688e-14  & &25600 & 148.8 \\ 
\hline
\multicolumn{11}{|c|}{HBVM(2,1)} \\
\hline
$h$          & $e_y$ & rate & $e_H$ & --- & $e_K$ & rate & $e_M$ & rate & blend & time \\
\hline
1.00e-01 &  2.184e-01 & --- &  1.110e-15  & & 4.652e-04 & --- & 1.397e-03  & --- & 6025   & 54.1 \\ 
5.00e-02 &  5.741e-02 & 1.9 &  9.159e-16  & & 1.171e-04 & 2.0 & 3.527e-04  & 2.0    & 9389   & 87.0\\
2.50e-02 &  1.453e-02 &  2.0    &  1.110e-15 & &  2.931e-05 & 2.0 &  8.839e-05  & 2.0    & 14441 & 138.6 \\
1.25e-02 &  3.642e-03 &  2.0    &  1.471e-15 & & 7.331e-06  & 2.0 & 2.211e-05  & 2.0    & 25600  & 257.5 \\
\hline
\multicolumn{11}{|c|}{HBVM(2,2)} \\
\hline
$h$          & $e_y$ & rate & $e_H$ & rate & $e_K$ & --- & $e_M$ & --- & blend & time \\
\hline
1.00e-01 & 3.184e-04 & ---  & 6.515e-07   &  --- & 8.549e-14 &&  1.688e-14 & &5606 & 48.6 \\ 
5.00e-02 & 1.995e-05 & 4.0 & 4.087e-08 & 4.0 & 7.883e-15 &&  1.155e-14 & &10329 & 90.3\\
2.50e-02 & 1.247e-06 & 4.0 & 2.559e-09   & 4.0 & 7.883e-15  && 9.770e-15 & &17414 & 154.2\\
1.25e-02 & 7.790e-08 & 4.0 & 1.599e-10 &  4.0 & 8.216e-15  && 1.155e-14 & &28419 & 258.0 \\ 
\hline
\multicolumn{11}{|c|}{HBVM(4,2)} \\
\hline
$h$          & $e_y$ & rate & $e_H$ & --- & $e_K$ & rate & $e_M$ & rate & blend & time \\
\hline
1.00e-01 & 2.853e-04 & ---  &   7.494e-16 & & 2.668e-07 & --- & 6.788e-07 & --- & 5979 & 58.1 \\ 
5.00e-02 & 1.789e-05 & 4.0 &  9.159e-16 & & 1.662e-08 & 4.0 & 4.250e-08 & 4.0 & 10363 & 103.4 \\
2.50e-02 & 1.119e-06 & 4.0 &  1.110e-15  & &  1.038e-09 & 4.0  & 2.658e-09 & 4.0 & 17508 & 177.9 \\
1.25e-02 & 6.989e-08 & 4.0 &  1.638e-15  & & 6.488e-11  & 4.0 & 1.661e-10  & 4.0 & 28606 & 298.7 \\
\hline
\multicolumn{11}{|c|}{HBVM(20,16)} \\
\hline
$h$          & $e_y$ & ---  & $e_H$ & --- & $e_K$ & --- & $e_M$ & ---  & blend & time \\
\hline
1.0           & 1.011e-10 &  & 5.551e-16 & & 8.993e-15 & & 7.550e-15 &  & 3473 & 60.7\\
\hline
\hline
\end{tabular}
}
\end{table}
All methods are implemented by the same Matlab code by suitably choosing the values of $k$ and $s$.
Consequently, all comparisons are quite fair. 
For each method, we list the maximum solution error (in infinity norm), $e_y$, the maximum Hamiltonian error, $e_H$, the maximum momentum error, $e_K$, and the maximum error in the total mass, $e_M$; when appropriate, we also list the corresponding estimated rate of convergence. Moreover, we also report the total number of {\em blended iterations} (\ref{blend}), and the measured execution times (in sec). From the obtained results, one deduces that:
\begin{enumerate}
\item the methods exhibit the correct convergence order in the solution;
\item the symplectic methods conserve the total mass and the momentum (which are both quadratic invariants), whereas the Hamiltonian error decreases with the prescribed order;
\item the energy-conserving methods conserve the Hamiltonian, whereas the total mass and momentum errors decrease with the prescribed order;
\item the energy-conserving methods are slightly more expensive than the corresponding symplectic methods of the same order, even though they are slightly more accurate;
\item the total number of {\em blended iterations}, for the given values of the time-step $h$ and of $s$, are quite similar; 
\item the higher order methods are more effective than the lower order ones;
\item the use of HBVMs as spectral methods in time (i.e., HBVM(20,16), in the present case) provides the best computational performance, with a practical conservation of all invariants, a very small solution error, and with a comparably small execution time, due to the possibility of using quite large time-step ($h=1$). This, in turn, further confirms what observed in \cite{BIMR2018}.
\end{enumerate}

\section{Conclusions and future developments}\label{fine}
In this paper we have developed a spectrally accurate algorithm, in both space and time, for numerically solving the initial value problem for the Manakov systems. The efficient implementation of the method has been also studied, and some numerical tests confirm both the theoretical findings and the effectiveness of the proposed approach. This paper follows a series of papers focused on the numerical solution of Hamiltonian PDEs, which we plan to further investigate in the future. 

\subsection*{Acknowledgements} This paper has been finalized during the visit of the second author at the Academy of Mathematics and Systems Science, Chinese Academy of Sciences, in Beijing, whose economic support is acknowledged. The work of Yifa Tang is partially supported by the National Natural Science Foundation of China (Grant No. 11771438). Beibei Zhu is supported by the National Natural Science Foundation of China (Grant No. 11901564), China Postdoctoral Science Foundation (Grant No. 2018M641505), and the National Center for Mathematics and Interdisciplinary Sciences, CAS.

\end{document}